\documentclass[11pt]{article}
\usepackage[ansinew]{inputenc}
\usepackage{array}
\usepackage{color}
\usepackage{amsmath}
\usepackage{amsxtra}
\usepackage{amstext}
\usepackage{amssymb}
\usepackage{latexsym}
\usepackage{amsfonts,cite}
\usepackage{graphics,epstopdf}
\usepackage{epsfig}
\begin{document}

\sloppy
\newcommand{\proof}{{\it Proof~}}
\newtheorem{thm}{Theorem}[section]
\newtheorem{cor}[thm]{Corollary}
\newtheorem{lem}[thm]{Lemma}
\newtheorem{prop}[thm]{Proposition}
\newtheorem{eg}[thm]{Example}
\newtheorem{defn}[thm]{Definition}

\newtheorem{rem}[thm]{Remark}
\numberwithin{equation}{section}

\thispagestyle{empty}
\begin{center}
{\Large{\bf Determinant approach to the 2-Iterated $q$-Appell and \\mixed type $q$-Appell polynomials}}${}^\star$\\~~

{\bf Subuhi Khan${}^\star${\footnote{${}^{\star}$Corresponding author;E-mail:~subuhi2006@gmail.com (Subuhi Khan)\\ \hspace*{0.45cm}Second author ;E-mail:~mumtazrst@gmail.com (Mumtaz Riyasat)}} and Mumtaz Riyasat}\\
Department of Mathematics\\
Aligarh Muslim University\\
Aligarh, India\\
\end{center}

\begin{abstract}
\noindent
In this article, the 2-iterated $q$-Appell family is introduced. Certain 2-iterated $q$-Appell and mixed type $q$-special polynomials are considered as members of this family. The numbers related to these polynomials are obtained. The determinant definitions for the 2-iterated $q$-Appell family and for the 2-iterated and mixed type $q$-Appell polynomials are established. The graphs of some 2-iterated $q$-Appell and mixed type $q$-Appell polynomials are drawn for different values of indices
and the roots of these polynomials are also investigated for certain values of index $n$ by using Matlab. Finally, the approximate solutions of the real zeros of these polynomials are given.
\end{abstract}



\section{ Introduction and preliminaries}

\parindent=8mm
The area of $q$-calculus has in the last twenty years served as a bridge between mathematics and physics.
Recently, there is a significant increase of activities in the area of $q$-calculus due to its applications in various fields
such as mechanics, mathematics and physics. The definitions and notations of $q$-calculus reviewed here are taken from \cite{qMath}.\\

The $q$-analogue of the shifted factorial $(a)_n$ is defined by
\begin{equation}
(a;q)_0=1,(a;q)_n=\prod\limits_{m=0}^{n-1}(1-q^{m}a),~n \in \mathbb{N}.
\end{equation}

The $q$-analogues of  a complex number $a$ and of the factorial function are defined by
\begin{equation}
[a]_q=\frac{1-q^a}{1-q},~~q \in \mathbb{C}-\{1\};~~~ a \in \mathbb{C},
\end{equation}
\begin{equation}
[n]_{q}!=\prod_{m=1}^n [m]_q=[1]_q [2]_q\cdots [n]_q=\frac{(q;q)_n}{(1-q)^n},~q\neq 1;~n \in \mathbb{N},~~[0]_q!=1,~q \in \mathbb{C};~0<q<1.
\end{equation}

The Gauss $q$-binomial coefficient ${n \brack k}_q$  is defined by
\begin{equation}
{n \brack k}_q=\frac{[n]_{q}!}{[k]_{q}![n-k]_{q}!}=\frac{(q;q)_n}{(q;q)_k (q;q)_{n-k}},~~~k=0,1,\ldots,n.
\end{equation}

The $q$-exponential function is defined as:
\begin{equation}
e_q(x)=\sum\limits_{n=0}^\infty \frac{x^n}{[n]_{q}!},~0 <|q|<1.
\end{equation}

The $q$-derivative $D_q f$ of a function $f$ at a point $0 \neq z \in \mathbb{C}$ is defined as:
\begin{equation}
D_{q}f(z)=\frac{f(qz)-f(z)}{qz-z},~~~~0<|q|<1.
\end{equation}

The $q$-analogue of Taylor series expansion of an arbitrary function $f(z)$ for $0<q<1$ is defined as:
\begin{equation}
f(z)=\sum\limits_{n=0}^\infty \frac{{(1-q)}^n}{(q;q)_n}D_q^n f(a)(z-a)^n_q,
\end{equation}
where $D_q^n f(a)$ is the $n^{th}$ $q$-derivative of the function $f$ at point $a$.

The class of Appell polynomials is characterized completely by Appell \cite{Appell} in 1880. In 1954, Sharma and Chak \cite{SharmaChak} introduced a $q$-analogue for the family of Appell polynomials and called this sequence of polynomials as $q$-harmonic. Further, in 1967 Al-Salaam \cite{Qappl} introduced the family of $q$-Appell polynomials $\{A_{n,q}(x)\}_{n\geq 0}$ and studied some of its properties. The $n$-degree polynomials $A_{n,q}(x)$ are called $q$-Appell provided they satisfy the following $q$-differential equation:
\begin{equation}
D_{q,x}~\{A_{n,q}(x)\}=[n]_q ~A_{n-1,q}(x),~n=0,1,2,\ldots~;~q \in \mathbb{C};~0<q<1.
\end{equation}

The $q$-Appell polynomials $A_{n,q}(x)$ are also defined by means of the following generating function \cite{Qappl}:
\begin{equation}
A_q(t)e_q(xt)=\sum\limits_{n=0}^\infty A_{n,q}(x)\frac{t^n}{[n]_{q}!},~~~~0<q<1,
\end{equation}
where
\begin{equation}
A_{q}(t):=\sum\limits_{n=0}^\infty A_{n,q}\frac{t^n}{[n]_q!},~~A_{0,q}=1;~A_q(t)\neq0.
\end{equation}

It is to be noted that $A_q(t)$ is an analytic function at $t=0$ and
\begin{equation}
A_{n,q}:=A_{n,q}(0)
\end{equation}
are the $q$-Appell numbers. Also, there exists a sequence of numbers $\{A_{n,q}\}_{n\geq0}$, such that the sequence $A_{n,q}(x)$ satisfies the following relation \cite{Qappl}:
\begin{equation}
A_{n,q}(x)=A_{n,q}+{n \brack 1}_q A_{n-1,q} ~x+{n \brack 2}_q A_{n-2,q}~ x^2+
\cdots+A_{0,q}~ x^n,~n=0,1,2,\ldots~.
\end{equation}

Based on appropriate selection for the function $A_q(t)$, different members belonging to the family of $q$-Appell polynomials
can be obtained. These members are listed in Table 1.\\

\noindent
\textbf{Table~1.~Certain members belonging to the $q$-Appell family} \\
\\
{\tiny{
\begin{tabular}{lllll}
\hline
&&&&\\
{\bf S.}  & {\bf Name of the} & $A_q(t)$ & {\bf Generating function} &  {\bf Series definition}\\
{\bf No.}  &  {\bf $q$-special} &   &  &\\
&{\bf polynomials and }&&&\\
&{\bf related numbers }&&&\\
\hline
{\bf I.} & $q$-Bernoulli &   $\left(\frac{t}{e_q(t)-1}\right)$  & $\left(\frac{t}{e_q(t)-1}\right) e_q(xt)=\sum\limits_{
n=0}^\infty B_{n,q}(x)\frac{t^n}{[n]_q!}$  & $B_{n,q}(x)=\sum\limits_{k=0}^n {n \brack k}_q B_{k,q} x^{n-k}$\\
& polynomials &&$\left(\frac{t}{e_q(t)-1}\right)=\sum\limits_{
n=0}^\infty B_{n,q}\frac{t^n}{[n]_q!}$ &\\
& and numbers &&$B_{n,q}:=B_{n,q}(0)$&\\
& \cite{Ernst1,QB}&&&\\
\hline

{\bf II.} & $q$-Euler&   $\left(\frac{2}{e_q(t)+1}\right)$  & $\left(\frac{2}{e_q(t)+1}\right)e_q(xt)=\sum\limits_{
n=0}^\infty E_{n,q}(x)\frac{t^n}{[n]_q!}$  & $E_{n,q}(x)=\sum\limits_{k=0}^n {n \brack k}_q E_{k,q} x^{n-k}$\\
& polynomials &&$\left(\frac{2}{e_q(t)+1}\right)=\sum\limits_{
n=0}^\infty E_{n,q}\frac{t^n}{[n]_q!}$ &\\
& and numbers &&$E_{n,q}:=E_{n,q}(0)$&\\
&\cite{Ernst1,qBE}&&&\\

\hline
{\bf III.} &$q$-Genocchi&   $\left(\frac{2 t}{e_q(t)+1}\right)$  & $\left(\frac{2 t}{e_q(t)+1}\right) e_q(xt)=\sum\limits_{
n=0}^\infty G_{n,q}(x)\frac{t^n}{[n]_q!}$  & $G_{n,q}(x)=\sum\limits_{k=0}^n {n \brack k}_q G_{k,q} x^{n-k}$\\
& polynomials &&$\left(\frac{2 t}{e_q(t)+1}\right)=\sum\limits_{
n=0}^\infty G_{n,q}\frac{t^n}{[n]_q!}$ &\\
& and numbers &&$G_{n,q}:=G_{n,q}(0)$&\\
& \cite{q-UmbralAppl,qBE}&&&\\

\hline
\end{tabular}}}\\
\vspace{.25cm}

Al-Salaam showed that the class of all $q$-Appell polynomials is a
maximal commutative subgroup of the group of all polynomial sets, {\em i.e.} the class of all $q$-Appell sequences is closed under the operation of $q$-umbral composition of polynomial sequences.

If $A_{n,q}(x)=\sum\limits_{k=0}^n a_{n,k;q}~x^k$ and $B_{n,q}(x)=\sum\limits_{k=0}^n b_{n,k;q}~x^k$ are sequences of $q$-polynomials, then the
$q$-umbral composition of $A_{n,q}(x)$ with $B_{n,q}(x)$ is defined to be the sequence
\begin{equation}
(A_{n,q}\circ B_q)(x)=\sum\limits_{k=0}^n a_{n,k;q}~B_{k,q}(x)=\sum\limits_{0\leq k\leq l\leq n}a_{n,k;q}b_{k,l;q}x^l.
\end{equation}

Under this operation the set of all $q$-Appell sequences is an abelian group and it can be seen by considering the fact that every $q$-Appell sequence is of the form
\begin{equation}
A_{n,q}(x) = \left(\sum_{k=0}^\infty {a_{k,q} \over [k]_q!} D_{q}^k\right) x^n
\end{equation}
and that the umbral composition of $q$-Appell sequences corresponds to multiplication of these formal $q$-power series in the operator $D_q$.

In 1982, Srivastava \cite{SrivqAppl} gave several characterizations for the well known Appell polynomials and for their basic analogues: the $q$-Appell
polynomials. In 1985, Roman proposed an approach similar to the umbral approach under the area of nonclassical umbral calculus which is called $q$-umbral calculus \cite{Rom,Rom1}. The $q$-Appell polynomials are studied using determinantal and umbral approaches,
see for example \cite{q-Appell,q-UmbralAppl}. Obtaining determinant forms for the $q$-Appell polynomials and their members is an
important aspect of such study. The determinant forms can be helpful for computation purposes and can also be useful in finding the solutions of general linear interpolation problems.\\

In this article, a family of the 2-iterated $q$-Appell polynomials is introduced. The determinant forms and other properties of this family and also for certain members belonging to this family are established. In Section 2, the 2-iterated $q$-Appell polynomials are introduced by means of generating function and series definition. Certain members belonging to the 2-iterated $q$-Appell family and some mixed type $q$-special polynomials and related numbers are also considered. In Section 3, the determinant forms of the 2-iterated and mixed type $q$-Appell polynomials are derived. In Section, 4, the roots of the 2-iterated $q$-Appell and mixed type $q$-Appell polynomials are investigated and their graphs are drawn for suitable values of indices.

\section{2-iterated $q$-Appell polynomials}

The 2-iterated $q$-Appell polynomials are introduced by means of generating function and series definition. In order to introduce the 2I$q$AP, two different sets of $q$-Appell polynomials $A_{n,q}^{I}(x)$ and ${A}_{n,q}^{II}(x)$ are considered.

Thus, from definitions (1.9) and (1.10), it follows that
\begin{equation}
A_{q}^{I}(t)e_q(xt)=\sum\limits_{n=0}^\infty A_{n,q}^{I}(x)\frac{t^n}{[n]_{q}!},~0<q<1,
\end{equation}
where
\begin{equation}
A_{q}^{I}(t):=\sum\limits_{n=0}^\infty A_{n,q}^{I}\frac{t^n}{[n]_q!};~A_{n,q}^{I}:=A_{n,q}^{I}(0);~A_{0,q}^{I}=1;~A_q^{I}(t)\neq 0
\end{equation}
and
\begin{equation}
A_{q}^{II}(t)e_q(xt)=\sum\limits_{n=0}^\infty {A}_{n,q}^{II}(x)\frac{t^n}{[n]_{q}!},~0<q<1,
\end{equation}
where
\begin{equation}
A_{q}^{II}(t):=\sum\limits_{n=0}^\infty {A}_{n,q}^{II}\frac{t^n}{[n]_q!};~{A}_{n,q}^{II}:={A}_{n,q}^{II}(0);~{A}_{0,q}^{II}=1;~A_q^{II}(t)\neq 0,\end{equation}
respectively.\\

The generating function for the 2I$q$AP is derived by using a different approach based on
replacement techniques. For this, the following result is proved:

\begin{thm} The following generating function for the 2-iterated $q$-Appell polynomials holds true:
\begin{equation}
G_q(x,t):=A_{q}^{I}(t)A_{q}^{II}(t)e_q(xt)=\sum\limits_{n=0}^\infty A_{n,q}^{[2]}(x)\frac{t^n}{[n]_q!},~0<q<1.
\end{equation}
\end{thm}

\noindent
\begin{proof} Expanding the $q$-exponential function $e_q(xt)$ in the l.h.s. of equation (2.1) and then replacing the powers of $x$, {\em i.e.} $x^0,~x^1,~x^2,\ldots,x^n$ by the corresponding polynomials
${A}_{0,q}^{II}(x),~{A}_{1,q}^{II}(x)$, $\ldots,{A}_{n,q}^{II}(x)$ in l.h.s. and replacing $x$ by the polynomial $A_{1,q}^{II}(x)$ in the r.h.s. of the resultant equation, it follows that
\begin{equation}
A_{q}^{I}(t)\Big[1+{A}_{1,q}^{II}(x)\frac{t}{[1]_q!}+{A}_{2,q}^{II}(x)\frac{t^2}{[2]_q!}+\ldots+{A}_{n,q}^{II}(x)\frac{t^n}{[n]_q!}+\ldots\Big]=
\sum\limits_{n=0}^\infty A_{n,q}^{I}\{A_{1,q}^{II}(x)\}\frac{t^n}{[n]_q!}.
\end{equation}

Summing up the series in l.h.s. and then using equation (2.3) and denoting the resultant 2I$q$AP in the r.h.s. by $A_{n,q}^{[2]}(x)$, that is
\begin{equation}
A_{n,q}^{[2]}(x)=A_{n,q}^{I}\{A_{1,q}^{II}(x)\},
\end{equation}
assertion (2.5) is proved.\\
\end{proof}

\noindent
{\bf Remark~2.1.}~The generating function (2.5) for the 2I$q$AP $A_{n,q}^{[2]}(x)$ is derived by replacing the powers of $x$
by the polynomials ${A}_{n,q}^{II}(x)~(n=0,1,\ldots)$ in generating function (2.1) of the $q$-Appell polynomials $A_{n,q}^{I}(x)$. The replacements of the powers of $x$
by the polynomials ${A}_{n,q}^{I}(x)~(n=0,1,\ldots)$ in generating function (2.3) of the $q$-Appell polynomials ${A}_{n,q}^{II}(x)$, yields the
same generating function.\\

\noindent
{\bf Remark~2.2.}~It is important to remark that equation (2.7) is the operational correspondence between the 2I$q$AP $A_{n,q}^{[2]}(x)$ and $q$-Appell polynomials $A_{n,q}(x)$.\\

\begin{thm} The following series definition for the 2-iterated $q$-Appell polynomials $A_{n,q}^{[2]}(x)$ holds true:
\begin{equation}
A_{n,q}^{[2]}(x)=\sum\limits_{k=0}^n {n \brack k}_q A_{k,q}^{I}~{A}_{n-k,q}^{II}(x).
\end{equation}
\end{thm}

\noindent
\begin{proof} Using equations (2.2) and (2.3) in the l.h.s. of generating function
(2.5) and then using Cauchy-product rule in the l.h.s. of the resultant equation, such that
\begin{equation}
\sum\limits_{n=0}^\infty \sum\limits_{k=0}^n {n \brack k}_q A_{k,q}^{I}~{A}_{n-k,q}^{II}(x)\frac{t^n}{[n]_q!}=\sum\limits_{n=0}^\infty A_{n,q}^{[2]}(x)\frac{t^n}{[n]_q!}.
\end{equation}

Equating the coefficients of same powers of $t$ in both sides of equation (2.9), assertion (2.8) follows.\\
\end{proof}

\noindent
{\bf Remark~2.3.}~In view of equations (1.14) and (1.15), it is remarked that the 2-iterated $q$-Appell polynomials $A_{n,q}^{[2]}(x)$ satisfy the following relation:
\begin{equation}
A_{n,q}^{[2]}(x)=\left(\sum\limits_{k=0}^\infty \frac{{A}_{k,q}^{I}}{[k]_q!}D_q^k\right)A_{n,q}^{II}(x).
\end{equation}

Certain members belonging to the $q$-Appell family are given in Table 1. Since, corresponding to each member belonging to the $q$-Appell family, there exists a new special polynomial belonging to the 2-iterated $q$-Appell family. Thus, by making suitable choice for the functions $A_q^{I}(t)$ and $A_q^{II}(t)$ in equations (2.5) and (2.8), the generating function and series definition for the corresponding member belonging to the 2-iterated $q$-Appell family can be obtained. These resultant 2-iterated $q$-Appell polynomials along with their notations, names, generating functions and series definitions
are given in Table 2.\\

\noindent
\textbf{Table~2.~Certain members belonging to the 2-iterated $q$-Appell family} \\
\\
{\tiny{
\begin{tabular}{lllll}
\hline
&&&&\\
{\bf S.}    & $A_q^{I}(t)=A_q^{II}(t)$ & {\bf Notation }  &  {\bf Generating function}  &  {\bf Series  definition}    \\
{\bf No.} &&{\bf and name}&&\\
 &  & {\bf of the  }&  &   \\
 &&{\bf  resultant}&&\\
&&{\bf 2I$q$AP}&&\\
\hline
{\bf I.} &$\left(\frac{t}{e_q(t)-1}\right)$   & ${B}_{n,q}^{[2]}(x)$:=&  $\left(\frac{t}{ e_q(t)-1}\right)^{2}e_q(xt)=\sum\limits_{n=0}^\infty {B}_{n,q}^{[2]}(x)\frac{t^n}{[n]_q!}$ &
${B}_{n,q}^{[2]}(x)=\sum\limits_{k=0}^n {n \brack k}_q {B}_{k,q}{B}_{n-k,q}(x)$\\

&& 2-iterated && $\hspace{2.5cm}$\\

&&$q$-Bernoulli&&\\
&& polynomials &&\\
&&(2I$q$BP) &&\\
\hline
{\bf II.} &$\left(\frac{2}{ e_q(t)+1}\right)$   & ${E}_{n,q}^{[2]}(x)$:=&  $\left(\frac{2}{ e_q(t)+1}\right)^{2}e_q(xt)=\sum\limits_{n=0}^\infty {E}_{n,q}^{[2]}(x)\frac{t^n}{[n]_q!}$ &
${E}_{n,q}^{[2]}(x)=\sum\limits_{k=0}^n {n \brack k}_q {E}_{k,q}{E}_{n-k,q}(x)$\\

&& 2-iterated && $\hspace{2.5cm}$\\

&&$q$-Euler&&\\
&& polynomials &&\\
&&(2I$q$EP) &&\\
\hline

{\bf III.} &$\left(\frac{2 t}{e_q(t)+1}\right)$   & ${G}_{n,q}^{[2]}(x)$:=&  $\left(\frac{2t}{ e_q(t)+1}\right)^{2}e_q(xt)=\sum\limits_{n=0}^\infty {G}_{n,q}^{[2]}(x)\frac{t^n}{[n]_q!}$ &
${G}_{n,q}^{[2]}(x)=\sum\limits_{k=0}^n {n \brack k}_q {G}_{k,q}{G}_{n-k,q}(x)$\\

&& 2-iterated && $\hspace{2.5cm}$\\
&&$q$-Genocchi&&\\
&& polynomials &&\\
&&(2I$q$GP) &&\\

\hline
\end{tabular}}}\\

\vspace{.35cm}

The combinations of any two different members of the $q$-Appell family in the 2-iterated $q$-Appell family, yields a new mixed type $q$-special
polynomial. Thus, by making suitable selections for the functions $A_q^{I}(t)$ and $A_q^{II}(t)$ in equations (2.5) and (2.8), the
generating function and series definition of the resultant mixed type $q$-special polynomials are obtained. The possible combinations of the $q$-Bernoulli, $q$-Euler and $q$-Genocchi polynomials (Table 1 (I-III)) are considered. The resultant mixed type $q$-special polynomials along with their notations, names, generating functions and series definitions are given in Table 3.\\

\noindent
\textbf{Table~3.~Certain mixed type $q$-special polynomials} \\
\\
{\tiny{
\begin{tabular}{lllll}
\hline
&&&&\\
{\bf S.}    & $A_q^{I}(t);~A_q^{II}(t)$ & {\bf Notation }  &  {\bf Generating functions}  &  {\bf Series  definitions}    \\
&&{\bf and} &&\\
{\bf No.} &&{\bf name}&&\\
 &  & {\bf of the  }&  &   \\
 &&{\bf mixed}&&\\
&&{\bf type }&&\\
&&{\bf $q$-special} &&\\
&&{\bf polynomials}&&\\
&&&&\\
\hline

{\bf I.} & $\left(\frac{t}{ e_q(t)-1}\right)$;  & $_B{E}_{n,q}(x):=$&  $\frac{2t}{( e_q(t)-1)( e_q(t)+1) }e_q(xt)=\sum\limits_{n=0}^\infty {}_B{E}_{n,q}(x)\frac{t^n}{[n]_q!}$ &
${}_B{E}_{n,q}(x)=\sum\limits_{k=0}^n {n \brack k}_q {E}_{k,q}{B}_{n-k,q}(x)$\\
&$\left(\frac{2}{ e_q(t)+1}\right)$&$q$-Bernoulli- &&\\

&&Euler&&\\
&& polynomials &&\\
&&($q$BEP) &&\\
\hline

{\bf II.} & $\left(\frac{t}{ e_q(t)-1}\right)$;  & $_B{G}_{n,q}(x)$:=&  $\frac{2t^2}{( e_q(t)-1) ( e_q(t)+1)}e_q(xt)=\sum\limits_{n=0}^\infty {}_B{G}_{n,q}(x)\frac{t^n}{[n]_q!}$ &
${}_B{G}_{n,q}(x)=\sum\limits_{k=0}^n {n \brack k}_q {G}_{k,q}{B}_{n-k,q}(x)$\\
&$\left(\frac{2 t}{e_q(t)+1}\right)$&$q$-Bernoulli- &&\\
&&Genocchi&&\\
&& polynomials &&\\
&&($q$BGP) &&\\

\hline
&&&&\\

{\bf III.} & $\left(\frac{2}{ e_q(t)+1}\right)$;  & $_E{G}_{n,q}(x)$:=& $\left(\frac{2 t^{1/2}}{ e_q(t)+1}\right)^{2}e_q(xt)=\sum\limits_{n=0}^\infty {}_E{G}_{n,q}(x)\frac{t^n}{[n]_q!}$  &
${}_E{G}_{n,q}(x)=\sum\limits_{k=0}^n {n \brack k}_q {G}_{k,q}{E}_{n-k,q}(x)$\\
&$\left(\frac{2 t}{ e_q(t)+1}\right)$&$q$-Euler-&&\\

&&Genocchi&&\\
&& polynomials &&\\
&&($q$EGP) &&\\
\hline
\end{tabular}}}\\
\vspace{.35cm}

\noindent
{\bf Remark~2.4.}
From Remark 2.1 and Table 3, it follows that
\begin{equation}
_{B}{E}_{n,q}(x)\equiv{}_E{B}_{n,q}(x);~~{}_B{G}_{n,q}(x)\equiv{}_G{B}_{n,q}(x);~~{}_E{G}_{n,q}(x)\equiv{}_G{E}_{n,q}(x),
\end{equation}
where ${}_E{B}_{n,q}(x)$, ${}_G{B}_{n,q}(x)$ and ${}_G{E}_{n,q}(x)$ are the $q$-Euler-Bernoulli polynomials ($q$EBP),
$q$-Genocchi-Bernoulli polynomials ($q$GBP) and $q$-Genocchi-Euler polynomials ($q$GEP). This is due to the fact that the set of all $q$-Appell sequences is closed under the operation of $q$-umbral composition of polynomial sequences and form an abelian group.\\

Next, the numbers related to the 2I$q$AP $A_{n,q}^{[2]}(x)$ and also related to certain members of this family given in Tables 2 and 3 are explored. As mentioned in definition (1.11), the $q$-Appell numbers are the values of the polynomials $A_{n,q}(x)$ at $x=0$
and satisfy equation (1.12).\\

Taking $x=0$ in equation (2.8), the following series definition for the 2-iterated $q$-Appell numbers (2I$q$AN) $A_{n,q}^{[2]}:=A_{n,q}^{[2]}(0)$ is obtained

\begin{equation}
A_{n,q}^{[2]}=\sum\limits_{k=0}^n {n \brack k}_q A_{k,q}^{I}A_{n-k,q}^{II},
\end{equation}
where $A_{n,q}^{I}$ and $A_{n,q}^{(2)}$ denote two different sets of $q$-Appell numbers.\\

Next, the numbers related to the 2I$q$BP $B_{n,q}^{[2]}(x)$, 2I$q$EP $E_{n,q}^{[2]}(x)$ and 2I$q$GP $G_{n,q}^{[2]}(x)$ are explored. Taking $x=0$ in series definitions of 2I$q$BP $B_{n,q}^{[2]}(x)$, 2I$q$EP $E_{n,q}^{[2]}(x)$ and 2I$q$GP $G_{n,q}^{[2]}(x)$
given in Table 2 (I-III) and in view of notations given in Table 1, the 2-iterated $q$-Bernoulli, $q$-Euler and $q$-Genocchi
numbers are obtained. These numbers are given in Table 4.\\

\noindent
{\bf Table 4.~2-iterated $q$-numbers} \\
\\
{\tiny{\begin{tabular}{lll}

         \hline
         &&\\
         {\bf S.No.} & {\bf Notation and name of the 2-iterated $q$-number} & {\bf Series definition} \\
         &&\\
         \hline

         {I.} & ${B}_{n,q}^{[2]}:={B}_{n,q}^{[2]}(0)$ & ${B}_{n,q}^{[2]}=\sum\limits_{k=0}^n {n \brack k}_q {B}_{k,q}{B}_{n-k,q}$ \\
         &2-iterated $q$-Bernoulli numbers (2I$q$BN)&\\
         \hline
         &&\\
         {II.} & ${E}_{n,q}^{[2]}:={E}_{n,q}^{[2]}(0)$& ${E}_{n,q}^{[2]}=\sum\limits_{k=0}^n {n \brack k}_q {E}_{k,q}{E}_{n-k,q}$ \\
         &2-iterated $q$-Euler numbers (2I$q$EN)&\\
         \hline
         &&\\
         {III.} & ${G}_{n,q}^{[2]}:={G}_{n,q}^{[2]}(0)$& ${G}_{n,q}^{[2]}=\sum\limits_{k=0}^n {n \brack k}_q {G}_{k,q}{G}_{n-k,q}$ \\
         &2-iterated $q$-Genocchi numbers (2I$q$GN) &\\
         \hline
         \end{tabular}}}\\
\vspace{.35cm}

Further, the numbers corresponding to the mixed type $q$-special polynomials given in Table 3 are explored. Taking $x=0$ in series definitions of the $q$BEP $_BE_{n,q}(x)$, $q$BGP $_BG_{n,q}(x)$ and $q$EGP $_EG_{n,q}(x)$ (Table 3 (I-III)) and
in view of notations given in Table 1, the $q$-Bernoulli-Euler, $q$-Bernoulli-Genocchi and $q$-Euler-Genocchi numbers are obtained.
These numbers are given in Table 5.\\

\newpage

\noindent
{\bf Table 5.~Mixed type $q$-numbers} \\
\\
{\tiny{
\begin{tabular}{lll}
\hline
&&\\
  {\bf S.No.} & {\bf Notation and name of the mixed type $q$-numbers} & {\bf Series definition} \\
  &&\\
  \hline
  {\bf I.} & ${}_B{E}_{n,q}:={}_B{E}_{n,q}(0)$  & ${}_B{E}_{n,q}=\sum\limits_{k=0}^n {n \brack k}_q {E}_{k,q}{B}_{n-k,q}$\\
&$q$-Bernoulli-Euler numbers ($q$BEN)&\\
  \hline
  &&\\
  {\bf II.} & ${}_B{G}_{n,q}:={}_B{G}_{n,q}(0)$  & ${}_B{G}_{n,q}=\sum\limits_{k=0}^n {n \brack k}_q {G}_{k,q}{B}_{n-k,q}$\\
&$q$-Bernoulli-Genocchi numbers ($q$BGN) &\\
  \hline
  &&\\
  {\bf III.} & ${}_E{G}_{n,q}:={}_E{G}_{n,q}(0)$& ${}_E{G}_{n,q}=\sum\limits_{k=0}^n {n \brack k}_q {G}_{k,q}{E}_{n-k,q}$\\
&$q$-Euler-Genocchi numbers ($q$EGN)  &\\
  \hline
\end{tabular}}}\\
\vspace{.35cm}

\noindent
{\bf Note.}~From Remark 2.4 and Table 5, it follows that
\begin{equation}
_{B}{E}_{n,q}\equiv{}_E{B}_{n,q}:={}_E{B}_{n,q}(0);~~{}_B{G}_{n,q}\equiv{}_G{B}_{n,q}:={}_G{B}_{n,q}(0);~~{}_E{G}_{n,q}\equiv{}_G{E}_{n,q}:={}_G{E}_{n,q}(0),
\end{equation}
where ${}_E{B}_{n,q}$, ${}_G{B}_{n,q}$ and ${}_G{E}_{n,q}$ are the $q$-Euler-Bernoulli numbers ($q$EBN),
$q$-Genocchi-Bernoulli numbers ($q$GBN) and $q$-Genocchi-Euler numbers ($q$GEN), respectively.\\

In the next section, the determinant forms for the 2-iterated $q$-Appell polynomials and mixed type $q$-special polynomials are established.

\section{Determinant approach}

Keleshteri and Mahmudov \cite{q-Appell} studied the $q$-Appell polynomials from determinant point of view. The determinant forms of the special polynomials are important for the computational and applied purposes. This fact provides motivation to establish the determinant definitions of the $q$-special polynomials introduced in previous section. In order to define the 2I$q$AP $A_{n,q}^{[2]}(x)$ by means of determinant, the following result is proved:

\begin{thm}
The following determinant form for the 2-iterated $q$-Appell polynomials $A_{n,q}^{[2]}(x)$ of degree $n$ holds true:
\begin{equation}
\begin{array}{l}
A_{0,q}^{[2]}(x)=\frac{1}{\beta_{0,q}},
\\
A_{n,q}^{[2]}(x)=\frac{{(-1)}^n}{{(\beta_{0,q})}^{n+1}}\left|\begin{array}{ccccccc}
 1  &  {A}_{1,q}^{II}(x)  &  {A}_{2,q}^{II}(x)  & \cdots &  {A}_{n-1,q}^{II}(x)  & {A}_{n,q}^{II}(x)\\
\\
 {\beta}_{0,q}  &  {\beta}_{1,q}  &  {\beta}_{2,q}  & \cdots &  {\beta}_{n-1,q}  &  {\beta}_{n,q}\\
\\
 0  &  {\beta}_{0,q}  &  {2  \brack  1}_q~{\beta}_{1,q}  & \cdots & {n-1 \brack 1}_q~ {\beta}_{n-2,q} & {n \brack 1}_q~{\beta}_{n-1,q}\\
\\
 0 & 0 & {\beta}_{0,q} & \cdots & {n-1 \brack 2}_q~{\beta}_{n-3,q} &  {n \brack 2}_q~{\beta}_{n-2,q}\\
. & . & . & \cdots & . & . \\
 . & . & . & \cdots & . & . \\
 0 & 0 & 0 & \cdots & {\beta}_{0,q} & {n \brack n-1}_q~{\beta}_{1,q}
\end{array} \right|,\end{array}
\end{equation}
where $n=1,2,\ldots~~$ and ${A}_{n,q}^{II}(x)~(n=0,1,2,\ldots)$ are the $q$-Appell polynomials of degree $n$; $\beta_{0,q}\neq 0$ and
\begin{equation}\begin{array}{l}
\beta_{0,q}=\frac{1}{A_{0,q}^{I}},\\
\beta_{n,q}=-\frac{1}{A_{0,q}^{I}}\Big(\sum\limits_{k=1}^n {n \brack k}_qA_{k,q}^{I}~\beta_{n-k,q}\Big),~~n=1,2,\ldots~.
\end{array}\end{equation}
\end{thm}

\noindent
\begin{proof}
Let $A_{n,q}^{[2]}(x)$ be a sequence of the 2I$q$AP defined by equation (2.5) and $A_{n,q}^{I}, ~\beta_{n,q}$, be two numerical sequences (the coefficients of $q$-Taylor's series expansions of functions) such that
\begin{equation}
A_{q}^{I}(t)=A_{0,q}^{I}+\frac{t}{[1]_{q}!}A_{1,q}^{I}+\frac{t^2}{[2]_{q}!}A_{2,q}^{I}+\cdots+\frac{t^n}{[n]_{q}!}A_{n,q}^{I}+\cdots,~
n=0,1,\ldots;~A_{0,q}^{I}\neq0,
\end{equation}
\begin{equation}
\hat{A}_{q}^{I}(t)=\beta_{0,q}+\frac{t}{[1]_{q}!}\beta_{1,q}+\frac{t^2}{[2]_{q}!}\beta_{2,q}+\cdots+\frac{t^n}{[n]_{q}!}\beta_{n,q}+\cdots,~
n=0,1,\ldots;~\beta_{0,q}\neq0,
\end{equation}
satisfying
\begin{equation}
A_{q}^{I}(t)\hat{A}_{q}^{I}(t)=1.
\end{equation}

Then, according to the Cauchy-product rule, it follows that
\begin{equation}
A_{q}^{I}(t)\hat{A}_{q}^{I}(t)=\sum\limits_{n=0}^\infty \sum\limits_{k=0}^n {n \brack k}_qA_{k,q}^{I}~\beta_{n-k,q}\frac{t^n}{[n]_q!},\nonumber
\end{equation}
 which gives
\begin{equation}
\sum\limits_{k=0}^n {n \brack k}_q A_{k,q}^{I}~\beta_{n-k,q}=\begin{cases}
1,~~~~~\textrm{if}~n=0,\cr
0,~~~~~\textrm{if}~n>0.\cr
\end{cases}
\end{equation}

Consequently, the following holds:
\begin{equation}
\begin{cases}
\beta_{0,q}=\frac{1}{A_{0,q}^{I}},\cr
\beta_{n,q}=-\frac{1}{A_{0,q}^{I}}\Big(\sum\limits_{k=1}^n {n \brack k}_q A_{k,q}^{I}~\beta_{n-k,q}\Big),~~n=1,2,\ldots~~.
\end{cases}
\end{equation}

Multiplication of equation (2.5) by $\hat{A}_{q}^{I}(t)$ gives
\begin{equation}
A_{q}^{I}(t)\hat{A}_{q}^{I}(t){A}_{q}^{II}(t)e_q(xt)=\hat{A}_{q}^{I}(t) \sum\limits_{n=0}^\infty A_{n,q}^{[2]}(x) \frac{t^n}{n!},\nonumber
\end{equation}
which in view of equations (3.4), (3.5) and (2.3) gives
\begin{equation}
\sum\limits_{n=0}^\infty {A}_{n,q}^{II}(x)\frac{ t^n}{[n]_{q}!}=\sum\limits_{n=0}^\infty A_{n,q}^{[2]}(x)\frac{t^n}{[n]_{q}!}
\sum\limits_{k=0}^\infty \beta_{k,q}\frac{t^k}{[k]_{q}!}.
\end{equation}

Again, multiplication of the series on the r.h.s. of equation (3.8) according to Cauchy-product rule leads to the following system of infinite
equations in the unknowns $A_{n,q}^{[2]}(x)~(n=0,1,\ldots$):

\begin{equation}
\begin{cases}
A_{0,q}^{[2]}(x)\beta_{0,q}=1,\cr
\\
A_{0,q}^{[2]}(x)\beta_{1,q}+A_{1,q}^{[2]}(x)\beta_{0,q}={A}_{1,q}^{II}(x),\cr
\\
A_{0,q}^{[2]}(x)\beta_{2,q}+{2 \brack 1}_q A_{1,q}^{[2]}(x)\beta_{1,q}+A_{2,q}^{[2]}(x)\beta_{0,q}={A}_{2,q}^{II}(x),\cr
\vdots\\
A_{0,q}^{[2]}(x)\beta_{n-1,q}+{n-1 \brack 1}_q A_{1,q}^{[2]}(x)\beta_{n-2,q}+\cdots+A_{n,q}^{[2]}(x)\beta_{0,q}={A}_{n-1,q}^{II}(x),\cr
\\
A_{0,q}^{[2]}(x)\beta_{n,q}+{n \brack 1}_q A_{1,q}^{[2]}(x)\beta_{n-1,q}+\cdots+A_{n,q}^{[2]}(x)\beta_{0,q}={A}_{n,q}^{II}(x),\cr
\vdots
\end{cases}
\end{equation}

First equation of system (3.9), proves the first part of assertion (3.1). Also, the special form of system (3.9) (lower triangular) allows to work out the unknowns $A_{n,q}^{[2]}(x)$. Operating with the first $n+1$ equations simply by applying the Cramer's rule, it follows that

\begin{equation}
\begin{array}{l}
A_{n,q}^{[2]}(x)=\frac{\left|\begin{array}{cccccc}
                                          \beta_{0,q} & 0 & 0 & \cdots & 0 & 1\\
                                          \\
                                          \beta_{1,q} & \beta_{0,q} & 0 & \cdots & 0 & {A}_{1,q}^{II}(x)\\
                                          \\
                                          \beta_{2,q} & {2 \brack 1}_q~\beta_{1,q}& \beta_{0,q} & \cdots & 0& {A}_{2,q}^{II}(x)\\
                                          \\
                                          . & . & . & \cdots & . & . \\

                                          \beta_{n-1,q} & {n-1 \brack 1}_q~\beta_{n-2,q} & {n-1 \brack 2}_q~\beta_{n-3,q} & \cdots & \beta_{0,q} &{A}_{n-1,q}^{II}(x)\\
                                          \\
                                          \beta_{n,q} & {n \brack 1}_q~\beta_{n-1,q} & {n \brack 2}_q~\beta_{n-2,q} & \cdots & {n \brack n-1}_q~\beta_{1,q} & {A}_{n,q}^{II}(x)\\
                                        \end{array}\right|}{\left|\begin{array}{cccccc}
           \beta_{0,q} & 0 & 0 & \cdots & 0 & 0\\
           \\
           \beta_{1,q} & \beta_{0,q} & 0 & \cdots & 0 & 0\\
           \\
           \beta_{2,q} &  {2 \brack 1}_q~\beta_{1,q} & \beta_{0,q} & \cdots & 0 & 0 \\

           . & . & . & \cdots & . & . \\
           \beta_{n-1,q} & {n-1 \brack 1}_q~\beta_{n-2,q} & {n-1 \brack 2}_q~\beta_{n-3,q} & \cdots & \beta_{0,q} & 0\\
           \\
           \beta_{n,q} & {n \brack 1}_q~\beta_{n-1,q} & {n \brack 2}_q~\beta_{n-2,q} & \cdots & {n \brack n-1}_q~\beta_{1,q} & \beta_{0,q}\\
         \end{array}\right|}\end{array},
         \end{equation}
where, $~n=1,2,\ldots~$.\\

Now, taking the transpose of the determinant in the numerator and expanding the determinant in the denominator, so that
\begin{equation}
\begin{array}{l}
A_{n,q}^{[2]}(x)=\frac{1}{(\beta_{0,q})^{n+1}}\left|\begin{array}{cccccc}
  \beta_{0,q} & \beta_{1,q} & \beta_{2,q} & \cdots & \beta_{n-1,q} & \beta_{n,q}\\
  \\
  0 & \beta_{0,q} & {2 \brack 1}_q~\beta_{1,q} & \cdots & {n-1 \brack 1}_q~\beta_{n-2,q} &  {n \brack 1}_q~\beta_{n-1,q}\\
\\
  0 & 0 & \beta_{0,q} & \cdots & {n-1 \brack 2}_q~\beta_{n-3,q} & {n \brack 2}_q~\beta_{n-2,q}\\
\\
. & . & . & \cdots & . & . \\
  0 & 0 & 0 & \cdots & \beta_{0,q} & {n \brack n-1}_q~\beta_{1,q}\\
\\
  1 & {A}_{1,q}^{II}(x) &  {A}_{2,q}^{II}(x)& \cdots & {A}_{n-1,q}^{II}(x) & {A}_{n,q}^{II}(x)
\end{array}\right|,\end{array}
\end{equation}
which after $n$ circular row exchanges, that is after moving the $i^{th}$ row to the $(i+1)^{th}$ position for $i=1,2,\ldots,n-1$, yields second part of assertion (3.1).\\
\end{proof}

\noindent
{\bf Remark~3.1.}~Since the polynomials $B_{n,q}(x)$, $E_{n,q}(x)$ and $G_{n,q}(x)$
are particular members of the $q$-Appell family $A_{n,q}(x)$. Therefore, the determinant forms of
$B_{n,q}(x)$, $E_{n,q}(x)$ and $G_{n,q}(x)$ can be obtained by giving suitable values to the coefficients $\beta_{0,q}$ and ${\beta}_{i,q}~(i=1,2,\cdots,n)$
in the definition of $A_{n,q}(x)$. \\

Taking $\beta_{0,q}=1$, ${\beta}_{i,q}=\frac{1}{[i+1]_q}~(i=1,2,\cdots,n)$ in the determinant definition
of the $q$-Appell polynomials $A_{n,q}(x)$ \cite[p.12(19)]{q-Appell}, the following determinant definition of the $q$-Bernoulli polynomials $B_{n,q}(x)$ is obtained:\\

\noindent
{\bf Definition~3.1.}
The $q$-Bernoulli polynomials $B_{n,q}(x)$ of degree $n$ are defined by
\begin{equation}
\begin{array}{l}
B_{0,q}(x)=1,\\
B_{n,q}(x)=(-1)^n\left|\begin{array}{ccccccc}
                       1 & x & x^2 & \cdots& x^{n-1} & x^n\\
                       \\
                       1 & \frac{1}{[2]_q} & \frac{1}{[3]_q} & \cdots& \frac{1}{[n]_q} & \frac{1}{[n+1]_q}\\
                       \\
                       0 & 1 & {2 \brack 1}_q~\frac{1}{[2]_q} & \cdots& {n-1 \brack 1}_q~\frac{1}{[n-1]_q} & {n \brack 1}_q~\frac{1}{[n]_q}\\
                       \\
                       0 & 0 & 1 & \cdots& {n-1 \brack 2}_q~\frac{1}{[n-2]_q} & {n \brack 2}_q~\frac{1}{[n-1]_q} \\
                       . & . & . & \cdots & . & . \\
                       . & . & . & \cdots & . & . \\
                       0 & 0 & 0 & \cdots & 1 & {n \brack n-1}_q~\frac{1}{[2]_q}
                     \end{array}\right|,~~n=1,2,\cdots~.\end{array}
                     \end{equation}

The particular case of definition (3.12), for $n=4$ is considered in \cite[p.250]{DetSheffer}.\\

Next, taking $\beta_{0,q}=1$, ${\beta}_{i,q}=\frac{1}{2}~(i=1,2,\cdots,n)$ in the determinant definition
of the $q$-Appell polynomials $A_{n,q}(x)$ \cite[p.12(19)]{q-Appell}, the following determinant definition of the $q$-Euler polynomials $E_{n,q}(x)$
is obtained:\\

\noindent
{\bf Definition~3.2.}The $q$-Euler polynomials $E_{n,q}(x)$ of degree $n$ are defined by
\begin{equation}\begin{array}{l}
E_{0,q}(x)=1,\\
E_{n,q}(x)=(-1)^n\left|\begin{array}{ccccccc}
                 1 & x & x^2 & \cdots& x^{n-1} & x^n\\
                 \\
                 1 & \frac{1}{2} & \frac{1}{2} & \cdots & \frac{1}{2} & \frac{1}{2}\\
                 \\
                 0 & 1 & \frac{1}{2}{2 \brack 1}_q & \cdots & \frac{1}{2}{n-1 \brack 1}_q & \frac{1}{2}{n \brack 1}_q\\
                 \\
                 0 & 0 & 1 &\cdots& \frac{1}{2}{n-1 \brack 2}_q & \frac{1}{2}{n \brack 2}_q \\
                 . & . & . & \cdots & . & .\\
                 . & . & . & \cdots & . & . \\
                 0 & 0 & 0 & \cdots & 1 & \frac{1}{2}{n \brack n-1}_q\\
               \end{array}\right|,~~n=1,2,\cdots~.\end{array}\end{equation}

Further, taking $\beta_{0,q}=1$, ${\beta}_{i,q}=\frac{1}{2[i+1]_q}~(i=1,2,\cdots,n)$ in the determinant definition
of the $q$-Appell polynomials $A_{n,q}(x)$ \cite[p.12(19)]{q-Appell}, the following determinant definition of the $q$-Genocchi polynomials $G_{n,q}(x)$
is obtained:\\

\noindent
{\bf Definition~3.3.}
The $q$-Genocchi polynomials $G_{n,q}(x)$ of degree $n$ are defined by
\begin{equation}
\begin{array}{l}
G_{0,q}(x)=1,\\
G_{n,q}(x)=(-1)^n\left|\begin{array}{ccccccc}
                       1 & x & x^2 & \cdots& x^{n-1} & x^n\\
                       \\
                       1 & \frac{1}{2[2]_q} & \frac{1}{2[3]_q} & \cdots& \frac{1}{2[n]_q} & \frac{1}{2[n+1]_q}\\
                       \\
                       0 & 1 & {2 \brack 1}_q~\frac{1}{2[2]_q} & \cdots& {n-1 \brack 1}_q~\frac{1}{2[n-1]_q} & {n \brack 1}_q~\frac{1}{2[n]_q}\\
                       \\
                       0 & 0 & 1 & \cdots& {n-1 \brack 2}_q~\frac{1}{2[n-2]_q} & {n \brack 2}_q~\frac{1}{2[n-1]_q} \\
                       . & . & . & \cdots & . & . \\
                       . & . & . & \cdots & . & . \\
                       0 & 0 & 0 & \cdots & 1 & {n \brack n-1}_q~\frac{1}{2[2]_q}
                     \end{array}\right|,~~n=1,2,\cdots~.\end{array}\end{equation}

\noindent
{\bf Remark~3.2.}
Taking suitable values of the coefficients $\beta_{0,q}$ and ${\beta}_{i,q}~(i=1,2,\cdots,n)$ in the determinant
definition of the 2I$q$AP family, the determinant definitions for the 2-iterated $q$-members belonging to this family and certain mixed type
$q$-special polynomials can be obtained.\\

First, the determinant definitions for the members of the 2I$q$AP given in Table 2 are obtained. Taking $\beta_{0,q}=1$, ${\beta}_{i,q}=\frac{1}{[i+1]_q}~(i=1,2,\cdots,n)$ and $A_{n,q}^{II}(x)=B_{n,q}(x)$ in equation (3.1), the  following determinant definition of the 2I$q$BP $B_{n,q}^{[2]}(x)$ is obtained:\\

\noindent
{\bf Definition~3.4.}
The 2-iterated $q$-Bernoulli polynomials $B_{n,q}^{[2]}(x)$ of degree $n$ are defined by
\begin{equation}
\begin{array}{l}
B_{0,q}^{[2]}(x)=1,\\
B_{n,q}^{[2]}(x)=(-1)^n\left|\begin{array}{ccccccc}
                       1  &  B_{1,q}(x)  &  B_{2,q}(x)  & \cdots &  B_{n-1,q}(x)  & B_{n,q}(x)\\
                       \\
                       1 & \frac{1}{[2]_q} & \frac{1}{[3]_q} & \cdots& \frac{1}{[n]_q} & \frac{1}{[n+1]_q}\\
                       \\
                       0 & 1 & {2 \brack 1}_q~\frac{1}{[2]_q} & \cdots& {n-1 \brack 1}_q~\frac{1}{[n-1]_q} & {n \brack 1}_q~\frac{1}{[n]_q}\\
                       \\
                       0 & 0 & 1 & \cdots& {n-1 \brack 2}_q~\frac{1}{[n-2]_q} & {n \brack 2}_q~\frac{1}{[n-1]_q} \\
                       . & . & . & \cdots & . & . \\
                       . & . & . & \cdots & . & . \\
                       0 & 0 & 0 & \cdots & 1 & {n \brack n-1}_q~\frac{1}{[2]_q}
                     \end{array}\right|,~~n=1,2,\cdots~,\end{array}\end{equation}
                     where $B_{n,q}(x)~(n=0,1,2,\ldots)$ are the $q$-Bernoulli polynomials of degree $n$.\\

Taking $\beta_{0,q}=1$, ${\beta}_{i,q}=\frac{1}{2}~(i=1,2,\cdots,n)$ and $A_{n,q}^{II}(x)=E_{n,q}(x)$ in equation (3.1), the following determinant definition of the 2I$q$EP $E_{n,q}^{[2]}(x)$ is obtained:\\

\noindent
{\bf Definition~3.5.}
The 2-iterated $q$-Euler polynomials $E_{n,q}^{[2]}(x)$ of degree $n$ are defined by
\begin{equation}
\begin{array}{l}
E_{0,q}^{[2]}(x)=1,\\
E_{n,q}^{[2]}(x)={(-1)}^n\left|\begin{array}{ccccccc}
                       1  &  E_{1,q}(x)  &  E_{2,q}(x)  & \cdots &  E_{n-1,q}(x)  & E_{n,q}(x)\\
                       \\
                       1 & \frac{1}{2} & \frac{1}{2} & \cdots & \frac{1}{2} & \frac{1}{2}\\
                       \\
                        0 & 1 & \frac{1}{2}{2 \brack 1}_q & \cdots & \frac{1}{2}{n-1 \brack 1}_q & \frac{1}{2}{n \brack 1}_q\\
                        \\
                        0 & 0 & 1 &\cdots& \frac{1}{2}{n-1 \brack 2}_q & \frac{1}{2}{n \brack 2}_q \\
                        . & . & . & \cdots & . & .\\
                        . & . & . & \cdots & . & . \\
                         0 & 0 & 0 & \cdots & 1 & \frac{1}{2}{n \brack n-1}_q\\
                         \end{array}\right|,~n=1,2,\cdots,\end{array}\end{equation}
                          where $E_{n,q}(x)~(n=0,1,2,\ldots)$ are the $q$-Euler polynomials of degree $n$.\\

Taking $\beta_{0,q}=1$, ${\beta}_{i,q}=\frac{1}{2[i+1]_q}~(i=1,2,\cdots,n)$ and $A_{n,q}^{II}(x)=G_{n,q}(x)$ in equation (3.1), the following determinant definition of the 2I$q$GP $G_{n,q}^{[2]}(x)$ is obtained:\\

\noindent
{\bf Definition~3.6.}
The 2-iterated $q$-Genocchi polynomials $G_{n,q}^{[2]}(x)$ of degree $n$ are defined by
\begin{equation}
\begin{array}{l}
G_{0,q}^{[2]}(x)=1,\\
G_{n,q}^{[2]}(x)=(-1)^n\left|\begin{array}{ccccccc}
                       1  &  G_{1,q}(x)  &  G_{2,q}(x)  & \cdots &  G_{n-1,q}(x)  & G_{n,q}(x)\\
                       \\
                       1 & \frac{1}{2[2]_q} & \frac{1}{2[3]_q} & \cdots& \frac{1}{2[n]_q} & \frac{1}{2[n+1]_q}\\
                       \\
                       0 & 1 & {2 \brack 1}_q~\frac{1}{2[2]_q} & \cdots& {n-1 \brack 1}_q~\frac{1}{2[n-1]_q} & {n \brack 1}_q~\frac{1}{2[n]_q}\\
                       \\
                       0 & 0 & 1 & \cdots& {n-1 \brack 2}_q~\frac{1}{2[n-2]_q} & {n \brack 2}_q~\frac{1}{2[n-1]_q} \\
                       . & . & . & \cdots & . & . \\
                       . & . & . & \cdots & . & . \\
                       0 & 0 & 0 & \cdots & 1 & {n \brack n-1}_q~\frac{1}{2[2]_q}
                     \end{array}\right|,~~n=1,2,\cdots~,\end{array}\end{equation}
                     where $G_{n,q}(x)~(n=0,1,2,\ldots)$ are the $q$-Genocchi polynomials of degree $n$.\\

Next, the determinant definitions for the mixed type $q$-special polynomials given in Table 3 are obtained.\\

Taking $\beta_{0,q}=1$, ${\beta}_{i,q}=\frac{1}{2}~(i=1,2,\cdots,n)$ and $A_{n,q}^{II}(x)=B_{n,q}(x)$ in equation (3.1), the following determinant definition of the $q$BEP $_BE_{n,q}(x)$ is obtained:\\

\noindent
{\bf Definition~3.7.}
The $q$-Bernoulli-Euler polynomials $_BE_{n,q}(x)$ of degree $n$ are defined by
\begin{equation}
\begin{array}{l}
_BE_{0,q}(x)=1,\\
_BE_{n,q}(x)=(-1)^n\left|\begin{array}{ccccccc}
                       1  &  B_{1,q}(x)  &  B_{2,q}(x)  & \cdots &  B_{n-1,q}(x)  & B_{n,q}(x)\\
                       \\
                       1 & \frac{1}{2} & \frac{1}{2} & \cdots & \frac{1}{2} & \frac{1}{2}\\
                       \\
                       0 & 1 & \frac{1}{2}{2 \brack 1}_q & \cdots & \frac{1}{2}{n-1 \brack 1}_q & \frac{1}{2}{n \brack 1}_q\\
                       \\
                       0 & 0 & 1 &\cdots& \frac{1}{2}{n-1 \brack 2}_q & \frac{1}{2}{n \brack 2}_q \\
                       . & . & . & \cdots & . & .\\
                 . & . & . & \cdots & . & . \\
                 0 & 0 & 0 & \cdots & 1 & \frac{1}{2}{n \brack n-1}_q\\
                     \end{array}\right|,~n=1,2,\cdots~,\end{array}\end{equation}
                     where $B_{n,q}(x)~(n=0,1,2,\ldots)$ are the $q$-Bernoulli polynomials of degree $n$.\\

Taking $\beta_{0,q}=1$, ${\beta}_{i,q}=\frac{1}{2[i+1]_q}~(i=1,2,\cdots,n)$ and $A_{n,q}^{II}(x)=B_{n,q}(x)$ in equation (3.1), the following determinant definition of the $q$BGP $_BG_{n,q}(x)$ is obtained:\\

\noindent
{\bf Definition~3.8.}
The $q$-Bernoulli-Genocchi polynomials $_BG_{n,q}(x)$ of degree $n$ are defined by
\begin{equation}
\begin{array}{l}
_BG_{0,q}(x)=1,\\
_BG_{n,q}(x)=(-1)^n\left|\begin{array}{ccccccc}
                       1  &  B_{1,q}(x)  &  B_{2,q}(x)  & \cdots &  B_{n-1,q}(x)  & B_{n,q}(x)\\
                       \\
                       1 & \frac{1}{2[2]_q} & \frac{1}{2[3]_q} & \cdots& \frac{1}{2[n]_q} & \frac{1}{2[n+1]_q}\\
                       \\
                       0 & 1 & {2 \brack 1}_q~\frac{1}{2[2]_q} & \cdots& {n-1 \brack 1}_q~\frac{1}{2[n-1]_q} & {n \brack 1}_q~\frac{1}{2[n]_q}\\
                       \\
                       0 & 0 & 1 & \cdots& {n-1 \brack 2}_q~\frac{1}{2[n-2]_q} & {n \brack 2}_q~\frac{1}{2[n-1]_q} \\
                       . & . & . & \cdots & . & . \\
                       . & . & . & \cdots & . & . \\
                       0 & 0 & 0 & \cdots & 1 & {n \brack n-1}_q~\frac{1}{2[2]_q}
                     \end{array}\right|,~n=1,2,\cdots~,\end{array}\end{equation}
                     where $B_{n,q}(x)~(n=0,1,2,\ldots)$ are the $q$-Bernoulli polynomials of degree $n$.\\

Taking $\beta_{0,q}=1$, ${\beta}_{i,q}=\frac{1}{2[i+1]_q}~(i=1,2,\cdots,n)$ and $A_{n,q}^{II}(x)=E_{n,q}(x)$  in equation (3.1), the following determinant definition of the $q$EGP $_EG_{n,q}(x)$ is obtained:\\

\noindent
{\bf Definition~3.9.}
The $q$-Euler-Genocchi polynomials $_EG_{n,q}(x)$ of degree $n$ are defined by
\begin{equation}\begin{array}{l}
_EG_{0,q}(x)=1,\\
_EG_{n,q}(x)=(-1)^n\left|\begin{array}{ccccccc}
                  1  &  E_{1,q}(x)  &  E_{2,q}(x)  & \cdots &  E_{n-1,q}(x)  & E_{n,q}(x)\\
                  \\
                 1 & \frac{1}{2[2]_q} & \frac{1}{2[3]_q} & \cdots& \frac{1}{2[n]_q} & \frac{1}{2[n+1]_q}\\
                       \\
                       0 & 1 & {2 \brack 1}_q~\frac{1}{2[2]_q} & \cdots& {n-1\brack 1}_q~\frac{1}{2[n-1]_q} & {n \brack 1}_q~\frac{1}{2[n]_q}\\
                       \\
                       0 & 0 & 1 & \cdots& {n-1 \brack 2}_q~\frac{1}{2[n-2]_q} & {n \brack 2}_q~\frac{1}{2[n-1]_q} \\
                       . & . & . & \cdots & . & . \\
                       . & . & . & \cdots & . & . \\
                       0 & 0 & 0 & \cdots & 1 & {n \brack n-1}_q~\frac{1}{2[2]_q}
          \end{array}\right|,~~n=1,2,\cdots~,\end{array}\end{equation}
          where $E_{n,q}(x)~(n=0,1,2,\ldots)$ are the $q$-Euler polynomials of degree $n$.\\

\noindent
{\bf Remark~3.3.}~Taking $x=0$ in determinant definitions (3.12)-(3.14) of $B_{n,q}(x)$, $E_{n,q}(x)$ and $G_{n,q}(x)$ and on expanding the determinants with respect to first row and using suitable notations from Table 1 (I-III), the determinant definitions of the $B_{n,q}$, $E_{n,q}$ and $G_{n,q}$ can be obtained.\\

\noindent
{\bf Remark~3.4.}~Taking $x=0$ in determinant definitions (3.15)-(3.20) of $B_{n,q}^{[2]}(x)$, $E_{n,q}^{[2]}(x)$ and $G_{n,q}^{[2]}(x)$, $_BE_{n,q}(x)$, $_BG_{n,q}(x)$ and $_EG_{n,q}(x)$ and then using suitable notations from Tables 4 and 5 (I-III), the determinant definitions of the $B_{n,q}^{[2]}$, $E_{n,q}^{[2]}$, $G_{n,q}^{[2]}$, $_BE_{n,q}(x)$, $_BG_{n,q}$ and $_EG_{n,q}$ can be obtained.\\

To give the applications of the operational correspondence introduced in Section 2, the following
identities for the $q$-Appell polynomials are considered \cite{q-Appell}.\\
\begin{equation}
\begin{array}{llll}
A_{n,q}(x)&=\frac{1}{\beta_{0,q}}\Big(x^n-\sum\limits_{k=0}^{n-1}~{n \brack k}_q~\beta_{n-k,q}~A_{k,q}(x)\Big),~~~n=1,2,\cdots~~~,\\
x^n&=\sum\limits_{k=0}^{n}~{n \brack k}_q \beta_{n-k,q}~A_{k,q}(x),~~~n=1,2,\cdots~~~.\\
\end{array}
\end{equation}

Replacing the powers of $x$, {\em i.e.} $x^1$ and $x^n$ by the corresponding polynomials
${A}_{1,q}(x)$ and ${A}_{n,q}(x)$ in above equations and then using operational representation (2.7)
in the resultant equations, the following identities for the 2I$q$AP $A_{n,q}^{[2]}(x)$ are obtained:
\begin{equation}
\begin{array}{lll}
A_{n,q}^{[2]}(x)&=\frac{1}{\beta_{0,q}}\Big({A}_{n,q}(x)-\sum\limits_{k=0}^{n-1}~{n \brack k}_q~\beta_{n-k,q}~A_{k,q}^{[2]}(x)\Big),~~~n=1,2,\cdots~~~,\\
{A}_{n,q}(x)&=\sum\limits_{k=0}^{n}~{n \brack k}_q \beta_{n-k,q}~A_{k,q}^{[2]}(x),~~~n=1,2,\cdots~~~.\\
\end{array}
\end{equation}

The above examples illustrate that the operational correspondence established in this article can be applied to derive the results for the newly introduced
$q$-special polynomials given in Tables 2 and 3 from the results of the corresponding member belonging to the $q$-Appell family.\\

In the next section, the shapes of the 2I$q$BP $B_{n,q}^{[2]}(x)$, 2I$q$EP $E_{n,q}^{[2]}(x)$, 2I$q$GP $G_{n,q}^{[2]}(x)$, $q$BEP $_BE_{n,q}(x)$, $q$BGP $_BG_{n,q}(x)$ and $q$EGP $_EG_{n,q}(x)$ are displayed. The zeros of these polynomials are also computed by using Matlab. \\

\section{Graphical representation and roots}

There has been increasing interest in solving mathematical problems with the aid of computers. Numerical investigation of the roots of certain
$q$-polynomials are considered in \cite{Ryoo1}. Also, the shapes of a new class of $q$-Bernoulli polynomials are explored in \cite{Ryoo}. Here, the plots of the 2I$q$BP $B_{n,q}^{[2]}(x)$, 2I$q$EP $E_{n,q}^{[2]}(x)$, 2I$q$GP $G_{n,q}^{[2]}(x)$, $q$BEP $_BE_{n,q}(x)$, $q$BGP $_BG_{n,q}(x)$ and $q$EGP $_EG_{n,q}(x)$ are drawn for $n=1,2,3,4$ and $q=\frac{1}{2}$ $(0<q<1)$. \\

This shows the four plots combined into one for each of these polynomials. For this, the values of the first five $B_{n,q}$, $E_{n,q}$ and $G_{n,q}$ are required. The values of
first five $B_{n,q}$, $E_{n,q}$ and $G_{n,q}$  \cite{Ernst1,q-UmbralAppl} are listed in Table 6.\\

\noindent
{\bf Table 6.~Values of first five $B_{n,q}$, $E_{n,q}$ and $G_{n,q}$}\\
\\
{\tiny{
        \begin{tabular}{llllll}
         \hline
         &&&&&\\
         $n$ & 0 & 1 & 2 & 3 &4 \\
         \hline
         &&&&&\\
         $B_{n,q}$ & 1 & $-(1+q)^{-1}$ & $q^2([3]_q!)^{-1}$ & $(1-q)q^3([2]_q)^{-1}([4]_q)^{-1}$ & $q^4 (1-q^2-2q^3-q^4+q^6)([2]_q^2 [3]_q [5]_q)^{-1}$\\
         \hline
         &&&&&\\
         $E_{n,q}$ & 1 & $-\frac{1}{2}$ & $\frac{1}{4}(-1+q)$ & $\frac{1}{8}(-1+2q+2q^2-q^3)$ & $\frac{1}{16}(q-1)[3]_q!(q^2-4q+1)$\\
         \hline
         &&&&&\\
          $G_{n,q}$ & 1 & $\frac{q}{1+q}$ & $\frac{-(q^{3}+3q^{2}+4q+3)}{(1+q)(1+q+q^2)}$ & $\frac{1}{(1+q)^2}(2q^3+q^2)$ & $\frac{1}{1+q^2}(\frac{2q^3+q^2}{(1+q)^2}+\frac{q^3+3q^2+4q+3}{(1+q)(1+q+q^2)})-\frac{q}{q+1}-\frac{1}{[5]_q}-1$\\
         &&&&&\\
          \hline
       \end{tabular}}}\\
\vspace{.35cm}

The expressions of the first five $B_{n,q}(x)$, $E_{n,q}(x)$ and $G_{n,q}(x)$ are obtained by making use of the values of the first five $B_{n,q}$, $E_{n,q}$ and $G_{n,q}$ in the series definitions given in Table 1 (I-III). The expressions of first five $B_{n,q}(x)$, $E_{n,q}(x)$ and $G_{n,q}(x)$ are listed in Table 7.\\

\noindent
{\bf Table 7.~Expressions of first five $B_{n,q}(x)$, $E_{n,q}(x)$ and $G_{n,q}(x)$}\\
\\
{\tiny{
\begin{tabular}{llllll}
  \hline
  &&&&&\\
 $n$ & 0 & 1 & 2 & 3 &4\\
 \hline
 &&&&&\\
  $B_{n,q}(x)$ & $1$ & $x-\frac{1}{1+q}$ & $x^2-\frac{[2]_q}{1+q}x+\frac{q^2}{[3]_q[2]_q}$ & $x^3-\frac{[3]_q x^2}{1+q}+\frac{q^2 x}{[2]_q}+\frac{(1-q)q^3}{[2]_q [4]_q}$ & $x^4-\frac{[4]_q}{1+q}x^3+\frac{[4]_q q^2}{[2]_q^2}x^2+\frac{(1-q)q^3}{[2]_q}x$\\
  &&&&&$+q^4 (1-q^2-2q^3-q^4+q^6)([2]_q^2 [3]_q [5]_q)^{-1}$\\
 \hline
 &&&&&\\
  $E_{n,q}(x)$ & $1$ & $x-\frac{1}{2}$ & $x^2-\frac{[2]_q}{2}x+\frac{1}{4}(-1+q)$ & $x^3-\frac{[3]_q}{2}x^2+\frac{[3]_q}{4}(-1+q)x$& $x^4-\frac{[4]_q}{2}x^3+\frac{[4]_q[3]_q(q-1)}{4[2]_q}x^2+$\\
  &&&& $+\frac{1}{8}
  (-1+2q+2q^2-q^3)$& $\frac{[4]_q(-1+2q+2q^2-q^3)}{8}x+\frac{(q-1)[3]_q!(q^2-4q+1)}{16}$\\
  \hline
  &&&&&\\
  $G_{n,q}(x)$ & $1$ & $x+\frac{q}{1+q}$ & $x^2+qx-\frac{(q^{3}+3q^{2}+4q+3)}{(1+q)(1+q+q^{2})}$ & $x^3+\frac{[3]_q q}{(1+q)}x^2-\frac{(q^{3}+3q^{2}+4q+3)}{(1+q)}x$& $x^4+\frac{[4]_q q}{q+1}x^3-\frac{[4]-q[3]_q (q^3+3q^2+4q+3)}{[2]_q (1+q)(1+q+q^2)}x^2$ \\
  &&&&$+\frac{(2q^{3}+q^{2})}{(1+q)}$& $\frac{[4]_q(2q^3+q^4)}{(1+q)^2}x+\frac{1}{1+q^2}(\frac{2q^3+q^2}{(1+q)^2}+\frac{q^3+3q^2+4q+3}{(1+q)(1+q+q^2)})$\\
  &&&&&$-\frac{q}{q+1}-\frac{1}{[5]_q}-1$\\
  &&&&&\\
  \hline
\end{tabular}}}\\
\vspace{.5cm}

By making appropriate substitutions from Tables 6 and 7 in the series definitions given in Tables 2 and 3 (I-III), the expressions of the first five 2I$q$BP $B_{n,q}^{[2]}(x)$, 2I$q$EP $E_{n,q}^{[2]}(x)$, 2I$q$GP $G_{n,q}^{[2]}(x)$, $q$BEP $_BE_{n,q}(x)$, $q$BGP $_BG_{n,q}(x)$ and $q$EGP $_EG_{n,q}(x)$ for $q=\frac{1}{2}$ are obtained. These expressions are given in Table 8.\\

\noindent
{\bf Table~8.~Expressions of first five $B_{n,1/2}^{[2]}(x)$, $E_{n,1/2}^{[2]}(x)$, $G_{n,1/2}^{[2]}(x)$, $_BE_{n,1/2}(x)$, $_BG_{n,1/2}(x)$ \hspace*{2cm}and $_EG_{n,1/2}(x)$}\\
\\
{\tiny{
\begin{tabular}{llllll}
\hline
&&&&&\\
$n$ & 0 & 1 & 2 & 3 &4\\
\hline
$B_{n,1/2}^{[2]}(x)$& $1$ &$x-\frac{4}{3}$& $x^2-2x+\frac{6}{7}$&$x^3-\frac{7}{3}x^2+\frac{3}{2}x-\frac{8}{45}$&$x^4-\frac{5}{2}x^3+\frac{45}{24}x^2-\frac{8}{24}x-\frac{221}{7812}$\\
&&&&&\\
\hline
$E_{n,1/2}^{[2]}(x)$& $1$ &$x-1$&$x^2-\frac{3}{2}x-\frac{1}{16}$&$x^3-\frac{7}{4}x^2+\frac{7}{32}x+\frac{5}{16}$&$x^4-\frac{15}{8}x^3+\frac{35}{128}x^2+\frac{125}{1024}x+
\frac{71}{1024}$\\
&&&&&\\
\hline
$G_{n,1/2}^{[2]}(x)$& $1$ &$x+\frac{2}{3}$&$x^2+x-\frac{181}{42}$&$x^3+\frac{7}{6}x^2-\frac{181}{24}x-\frac{37}{18}$&$x^4+\frac{5}{4}x^3-\frac{905}{96}x^2-\frac{585}{144}x+
\frac{358499}{31248}$\\
&&&&&\\
\hline
$_BE_{n,1/2}(x)$ & $1$ &$x-\frac{7}{6}$&$x^2-\frac{3}{2}x+\frac{51}{168}$&$x^3-\frac{49}{24}x^2+\frac{79}{96}x+\frac{379}{2880}$&$x^4-\frac{35}{16}x^3-\frac{445}{384}x^2-
\frac{461}{1536}x-\frac{402305}{9999360}$\\
&&&&&\\
\hline
$_BG_{n,1/2}(x)$& $1$ &$x-\frac{1}{3}$&$x^2-\frac{1}{2}x-\frac{52}{21}$&$x^3-\frac{7}{12}x^2-\frac{5}{12}x+\frac{21}{70}$&$x^4-\frac{15}{24}x^3-\frac{5460}{1008}x^2+
\frac{5502}{1008}x-\frac{22203615}{27888840}$\\
&&&&&\\
\hline
$_EG_{n,1/2}(x)$&$1$ &$x-\frac{1}{6}$&$x^2-\frac{1}{4}x-\frac{439}{168}$&$x^3-\frac{7}{24}x^2-\frac{439}{96}x+\frac{1187}{576}$&$x^4-\frac{5}{16}x^3-\frac{595}{384}x^2+
\frac{3553}{32256}x+\frac{1226437}{19999872}$\\
&&&&&\\
\hline
\end{tabular}}}\\
\vspace{.25cm}

With the help of Matlab and by using the expressions for the first five $B_{n,1/2}^{[2]}(x)$, $E_{n,1/2}^{[2]}(x)$, $G_{n,1/2}^{[2]}(x)$, $_BE_{n,1/2}(x)$, $_BG_{n,1/2}(x)$ and $_EG_{n,1/2}(x)$ from Table 8 for $n=1,2,3,4$, the following graphs are drawn:\\

\begin{center}
\includegraphics[width=7.5cm,scale=1.5]{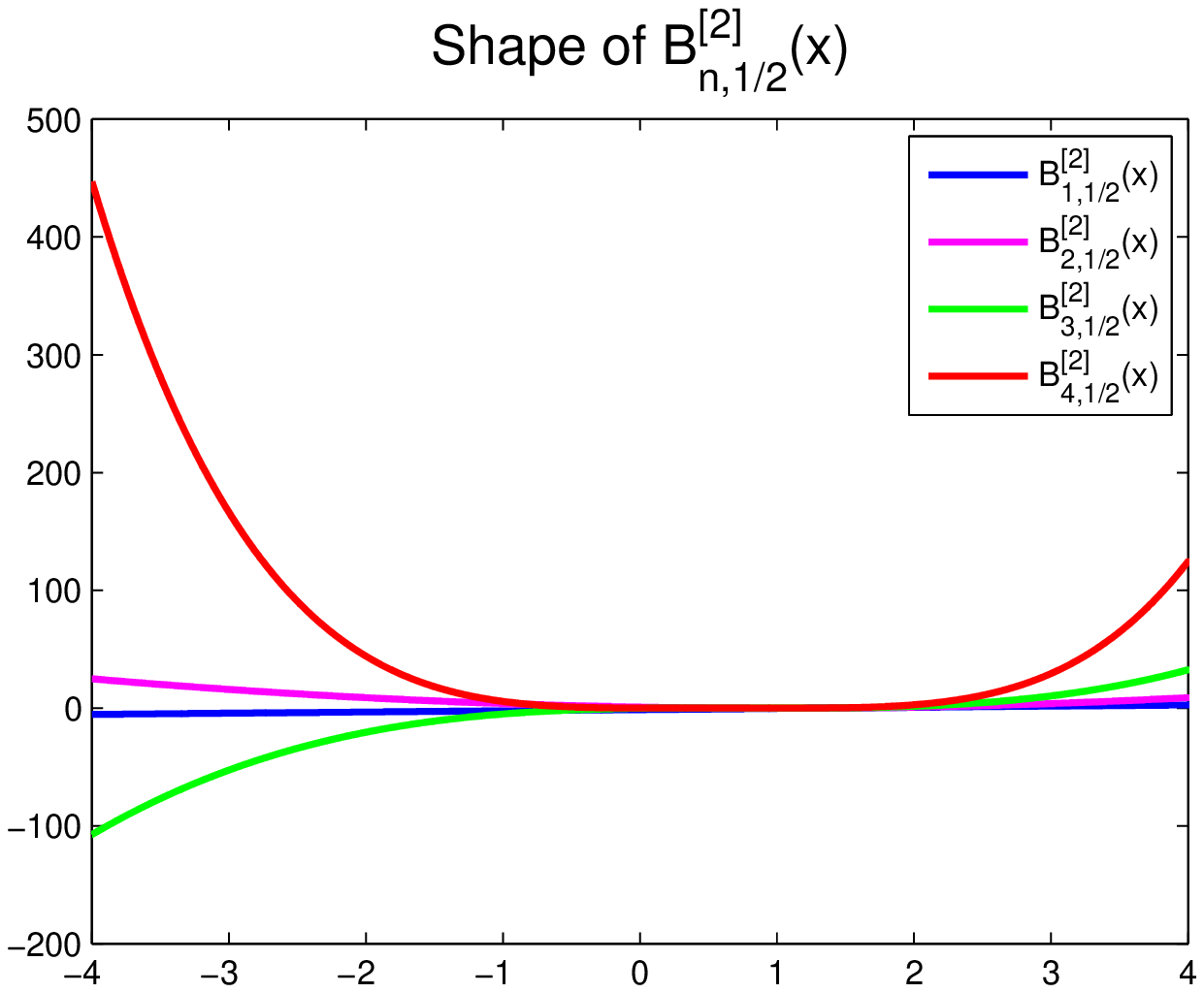}
\includegraphics[width=7.5cm,scale=1.5]{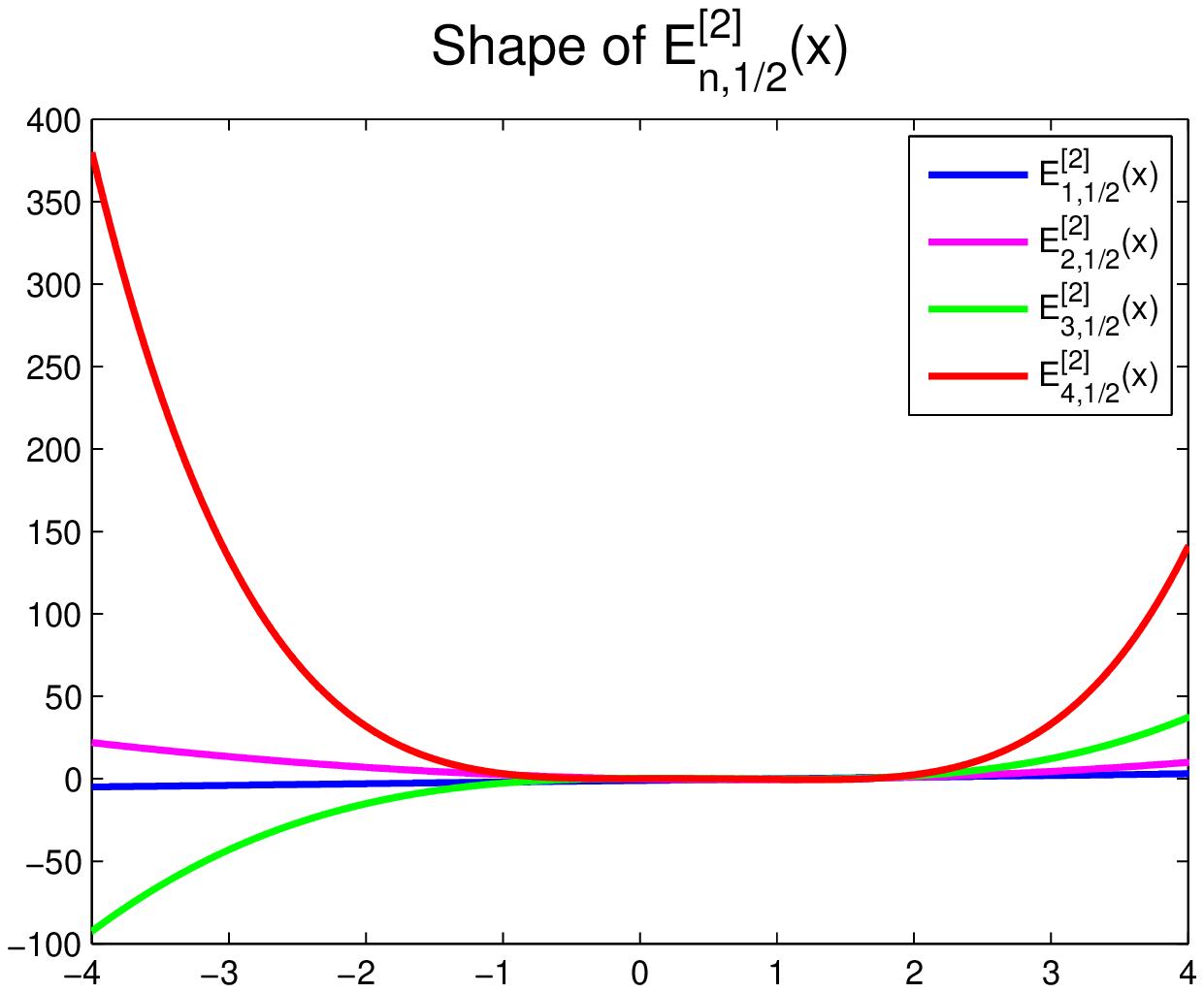}
{\bf Figure 4.1}\hspace{5cm}{\bf Figure 4.2}
\end{center}

\begin{center}
\includegraphics[width=7.5cm,scale=1.5]{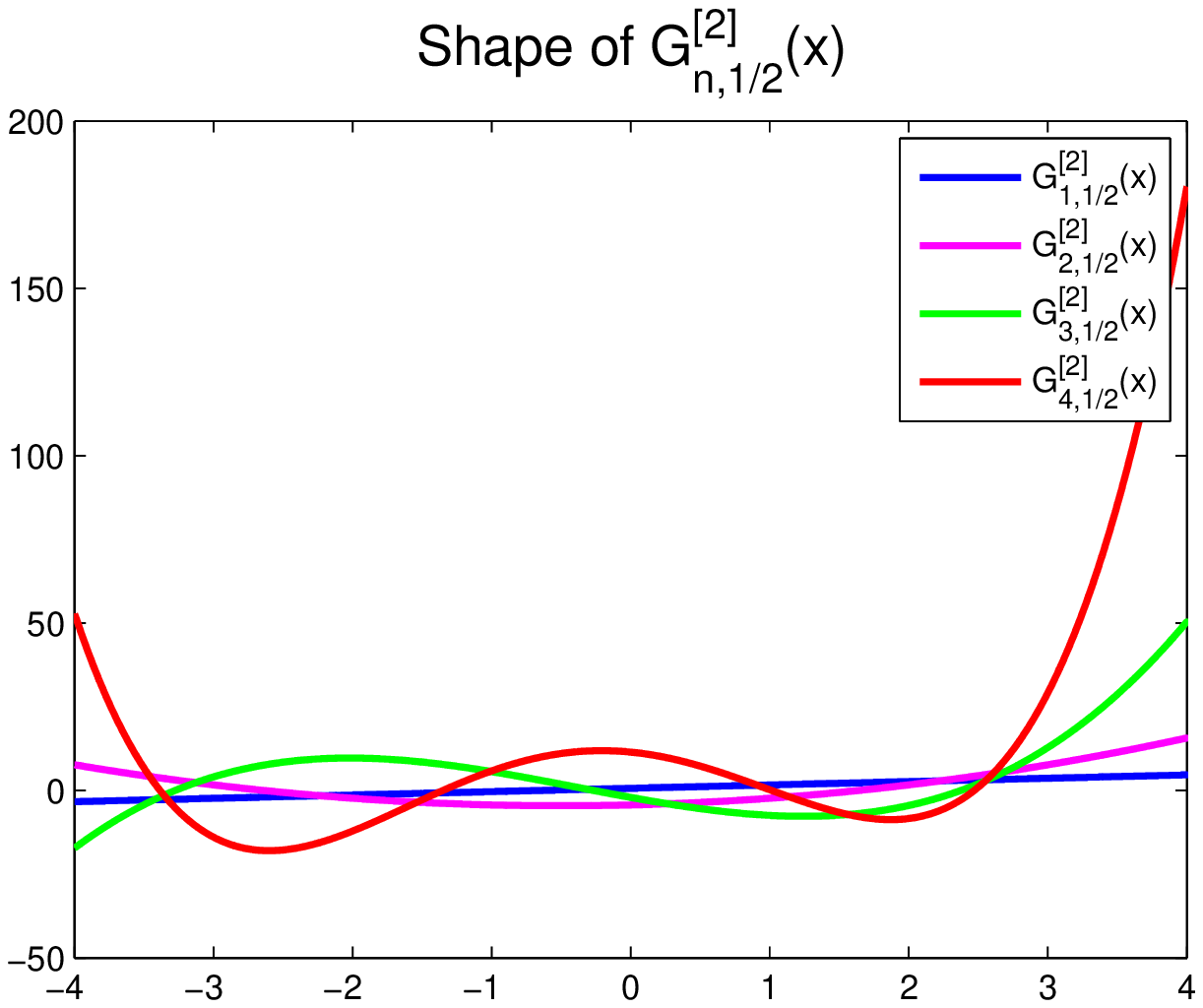}
\includegraphics[width=7.5cm,scale=1.5]{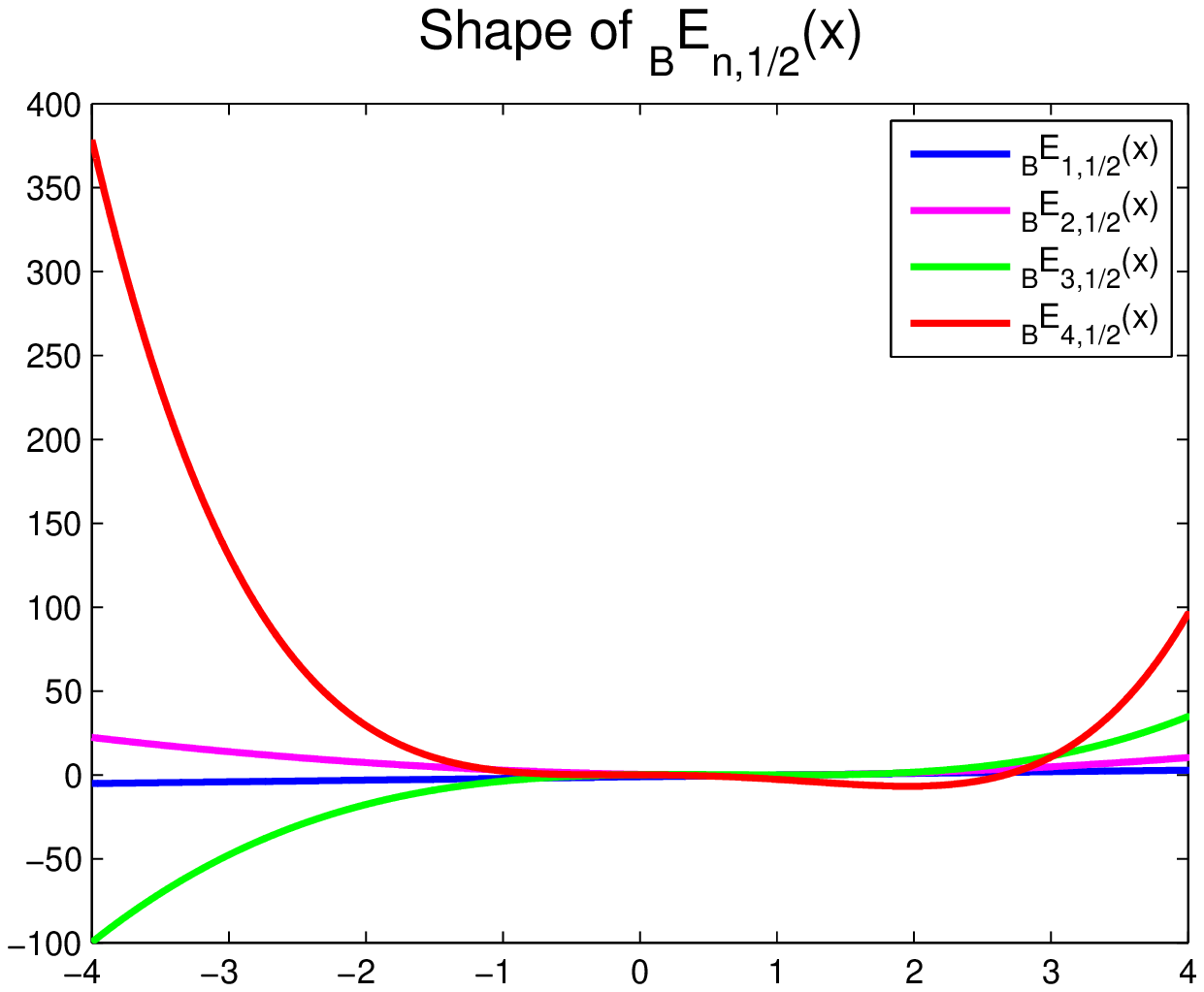}
{\bf Figure 4.3}\hspace{5cm}{\bf Figure 4.4}
\end{center}

\begin{center}
\includegraphics[width=7.5cm,scale=1.5]{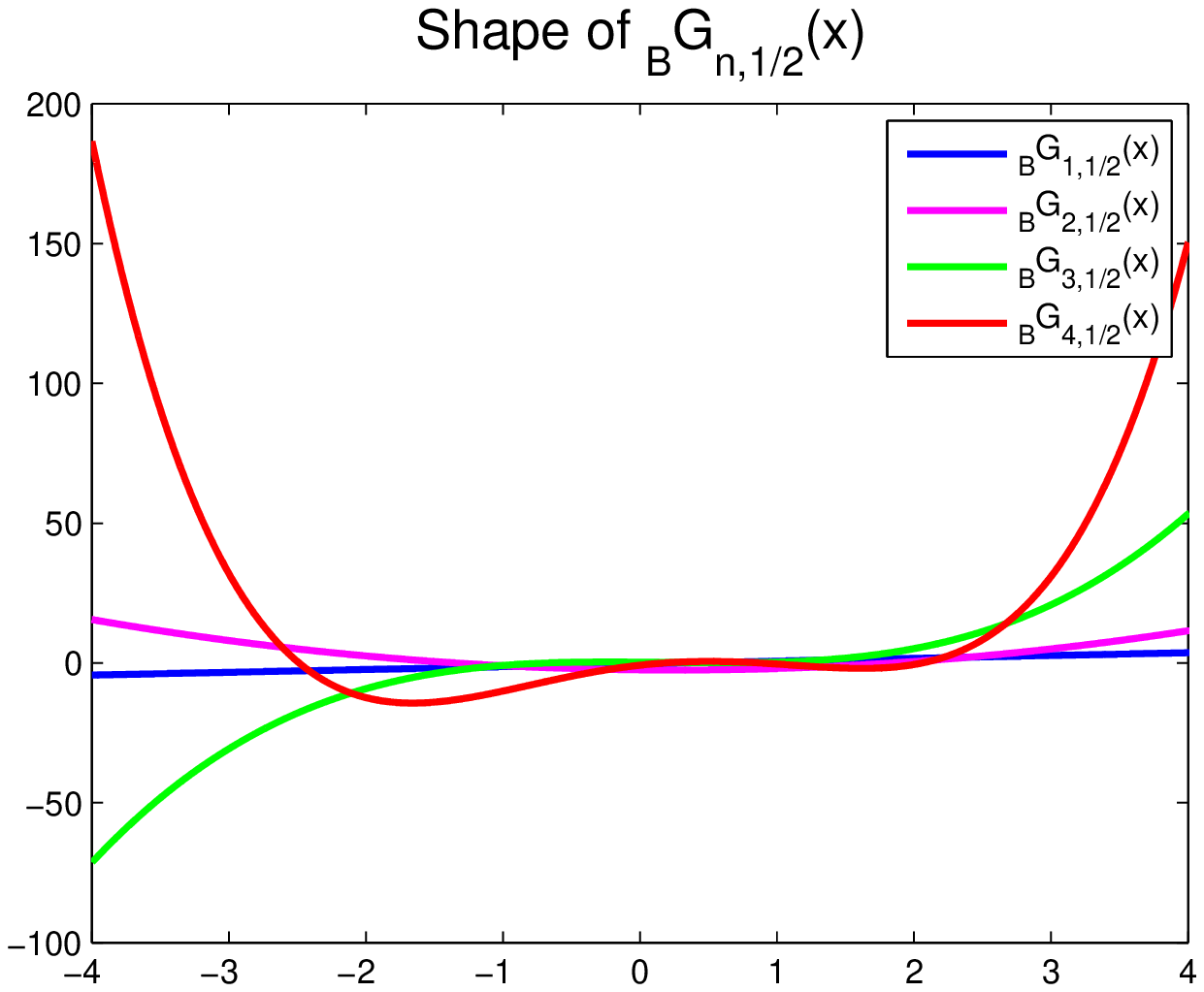}
\includegraphics[width=7.5cm,scale=1.5]{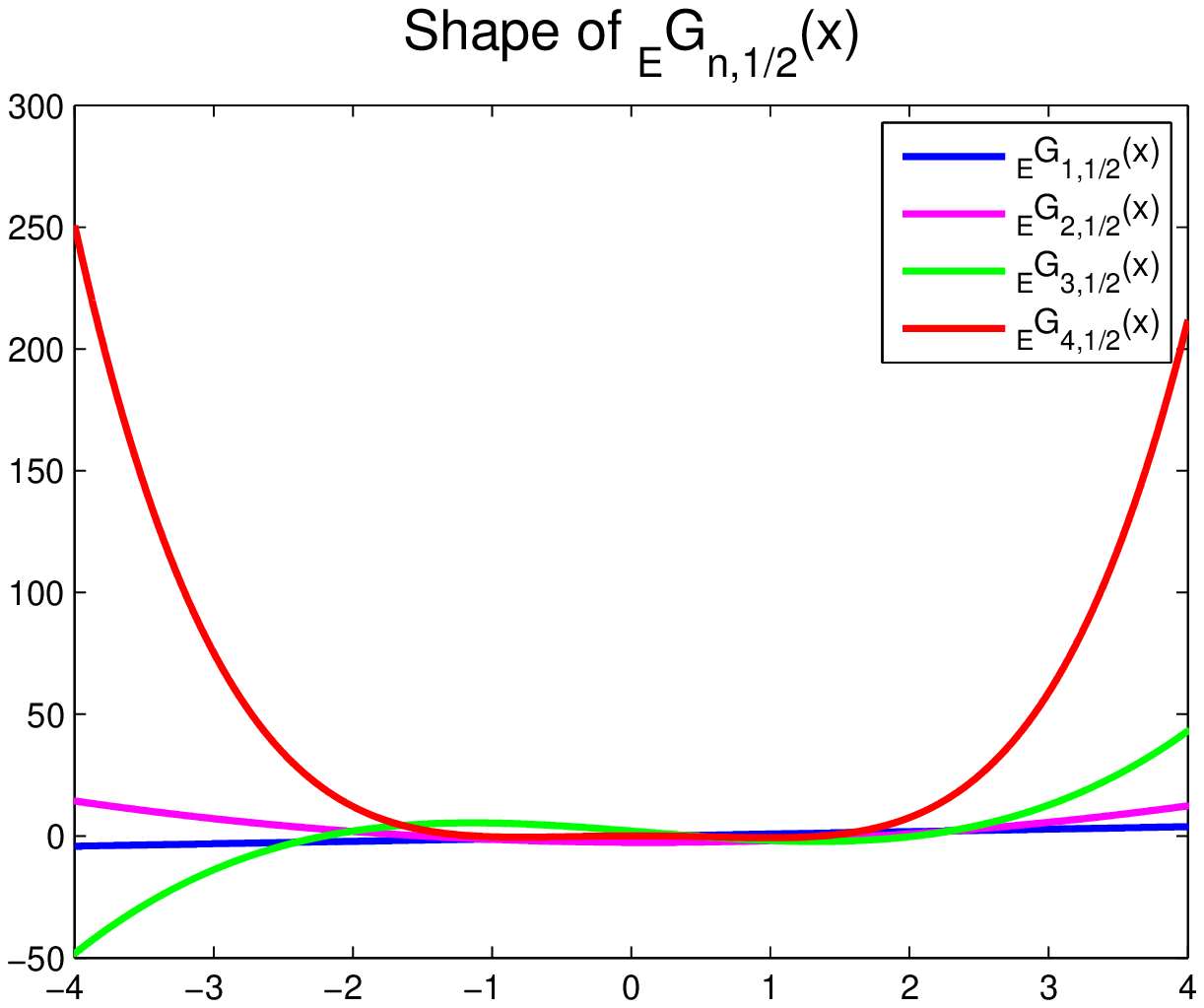}
{\bf Figure 4.5}\hspace{5cm}{\bf Figure 4.6}
\end{center}

Next, we find the real and complex zeros of the polynomials given in Table 8. The computation of these zeros
by hand is too complicated, therefore, we use Matlab to investigate these zeros. The investigation in this direction
will lead to a new approach employing numerical methods in the field of these $q$-special polynomials to appear in
mathematics and physics. \\

The real zeros of $B_{n,1/2}^{[2]}(x)$, $E_{n,1/2}^{[2]}(x)$, $G_{n,1/2}^{[2]}(x)$, $_BE_{n,1/2}(x)$, $_BG_{n,1/2}(x)$ and $_EG_{n,1/2}(x)$ are computed by using Matlab. These zeros are given in Table 9.\\

\noindent
{\bf Table 9.~Real zeros of $B_{n,1/2}^{[2]}(x)$, $E_{n,1/2}^{[2]}(x)$, $G_{n,1/2}^{[2]}(x)$, $_BE_{n,1/2}(x)$, $_BG_{n,1/2}(x)$ and $_EG_{n,1/2}(x)$}\\
\\
{\tiny{
\begin{tabular}{lllllll}
\hline
&&&&&&\\
$\textrm{Degree}~n$ & $B_{n,1/2}^{[2]}(x)$ & $E_{n,1/2}^{[2]}(x)$ & $G_{n,1/2}^{[2]}(x)$ & $_BE_{n,1/2}(x)$ & $_BG_{n,1/2}(x)$ & $_EG_{n,1/2}(x)$\\
&&&&&&\\
\hline
&&&&&&\\
$1$& $1.3333$ &$1.0000$ &$-0.6667$ &$1.1667$ &$0.3333$& $0.1667$\\
\hline
&&&&&&\\
$2$& $0.6220$, $1.3780$& $-0.0406$, $1.5406$ &$-2.6353$, $1.6353$ &$0.2411$, $1.2589$ &$-1.3433$, $1.8433$& $-1.4963$, $1.7463$\\
\hline
&&&&&&\\
$3$& $0.1522$, $0.9446$,  & $-3.2768$, $-0.2642$, & $-3.2768$, $-0.2642$,  & $-0.1213$, $0.7910$,  & $-0.6795$& $-2.2055$, $0.4583$, \\
&$1.2365$& $2.3743$ & $2.3743$  &$1.3719$ & & $2.0388$\\
\hline
&&&&&&\\
$4$& $-0.0617$, $0.3823$ &$0.5479$, $1.6488$ &$-3.3847$, $-1.3512$,  &$-0.2476$, $2.6664$ &$-2.4871$, $0.1776$, & $-1.1179$, $-0.1703$, \\
&&& $1.0514$, $2.4705$& &  $0.8751$, $2.0594$ &$0.2360$, $1.3647$\\
\hline
\end{tabular}}}\\
\vspace{.35cm}

Again, with the help of Matlab, we find the complex zeros of these polynomials. These complex zeros are given in Table 10.\\

\noindent
{\bf Table 10.~Complex zeros of $B_{n,1/2}^{[2]}(x)$, $E_{n,1/2}^{[2]}(x)$, $G_{n,1/2}^{[2]}(x)$, $_BE_{n,1/2}(x)$, $_BG_{n,1/2}(x)$ and \hspace*{2.0cm}$_EG_{n,1/2}(x)$}\\
\\
{\tiny{
\begin{tabular}{lllllll}
\hline
$\textrm{Degree}~n$ & $B_{n,1/2}^{[2]}(x)$ & $E_{n,1/2}^{[2]}(x)$ & $G_{n,1/2}^{[2]}(x)$ & $_BE_{n,1/2}(x)$ & $_BG_{n,1/2}(x)$ & $_EG_{n,1/2}(x)$\\
\hline
$1$&$-$& $-$& $-$& $-$& $-$& $-$\\
\hline

$2$&$-$& $-$& $-$& $-$& $-$& $-$\\
\hline

$3$& $-$& $-$&$-$ & $-$& -$0.6314+0.2068i$,& $-$\\
&&&&&$0.6314-0.2068i$&\\
\hline

$4$&$1.0897+0.1112i$,& $-0.1609+0.225i$,&$-$ & $-0.1157+0.2181i$,& $-$&$-$ \\
&$1.0897-0.1112i$&  $-0.1609-0.225i$  & &$-0.1157-0.2181i$& & \\
\hline
\end{tabular}}}\\
\vspace{.15cm}

Further, we plot the zeros of the polynomials $B_{n,1/2}^{[2]}(x)$, $E_{n,1/2}^{[2]}(x)$, $G_{n,1/2}^{[2]}(x)$, $_BE_{n,1/2}(x)$, $_BG_{n,1/2}(x)$ and $_EG_{n,1/2}(x)$ for $n=1,2,3,4$.\\

\begin{center}
\includegraphics[width=7.5cm,scale=1.5]{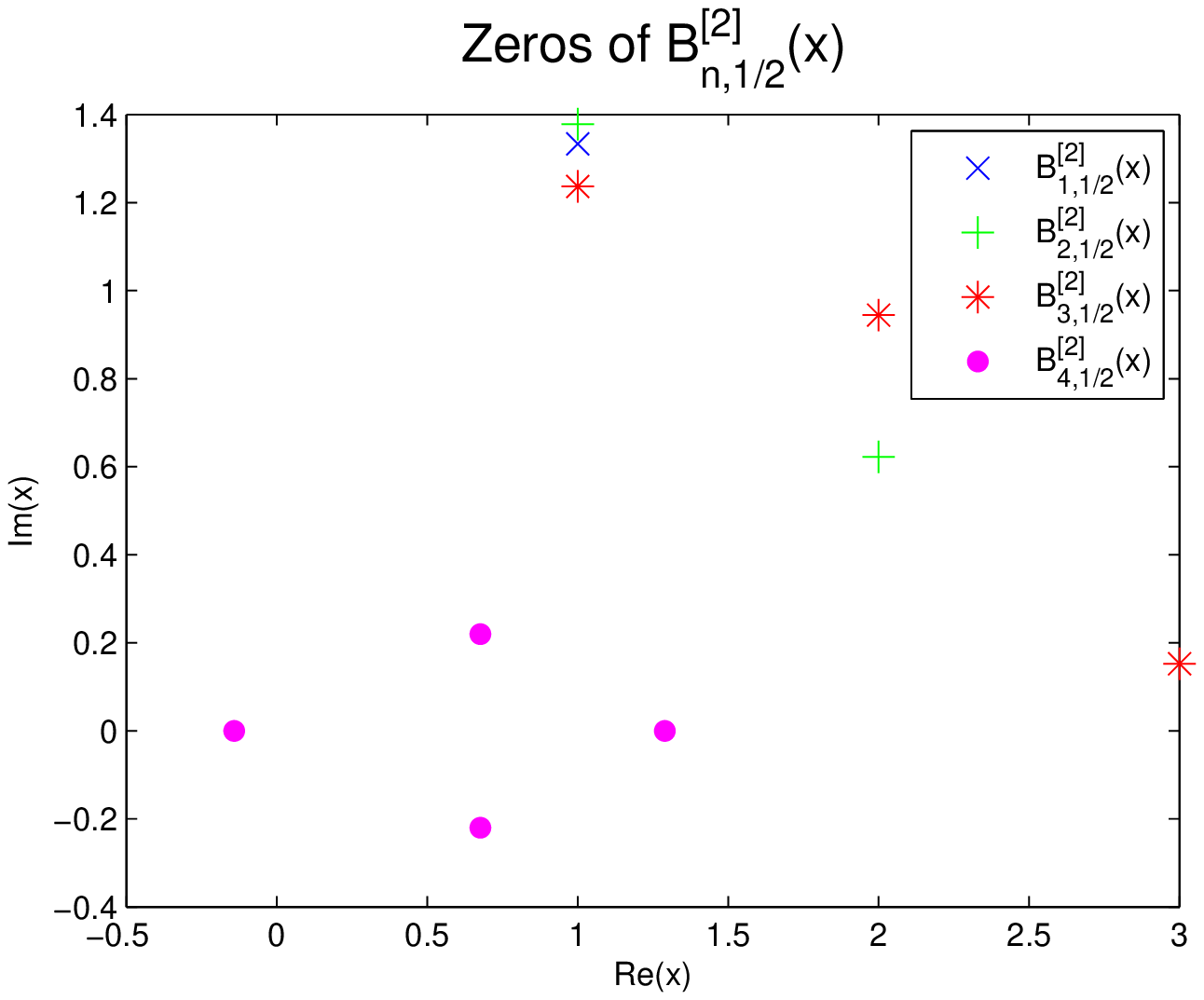}
\includegraphics[width=7.5cm,scale=1.5]{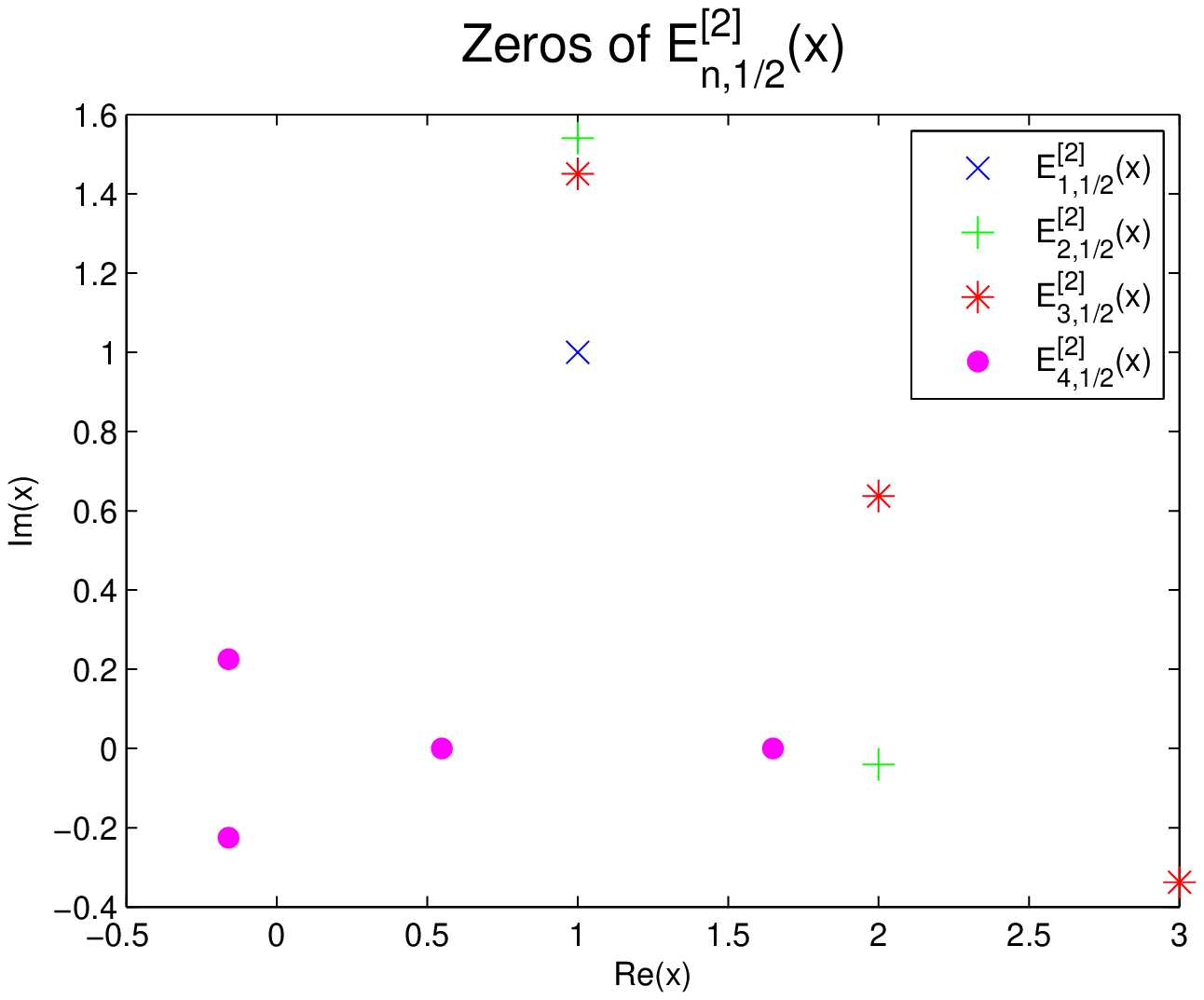}
{\bf Figure 4.7}\hspace{5cm}{\bf Figure 4.8}
\end{center}

\vspace{.35cm}
\begin{center}
\includegraphics[width=7.5cm,scale=1.5]{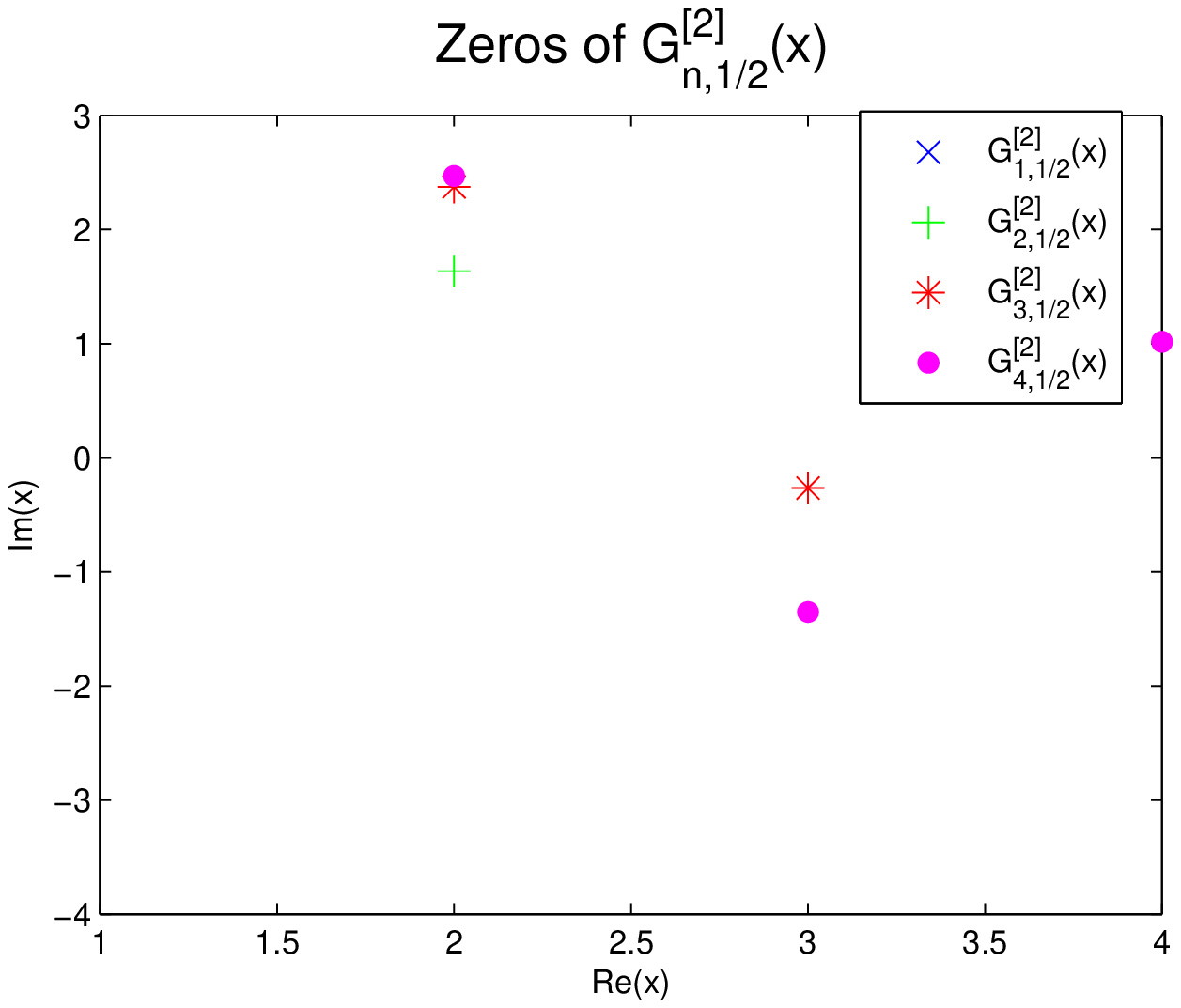}
\includegraphics[width=7.5cm,scale=1.5]{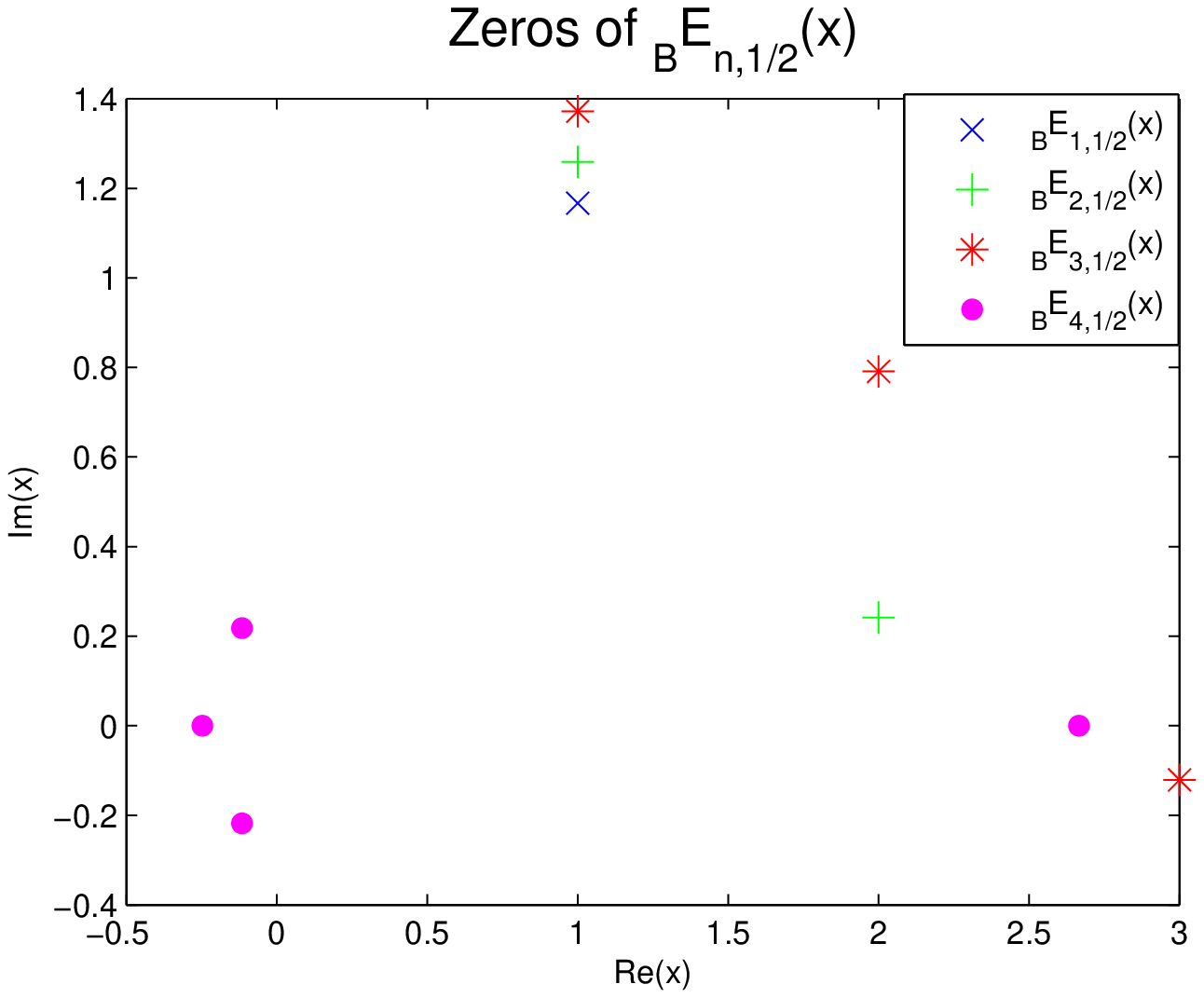}
{\bf Figure 4.9}\hspace{5cm}{\bf Figure 4.10}
\end{center}
\vspace{.35cm}

\begin{center}
\includegraphics[width=7.5cm,scale=1.5]{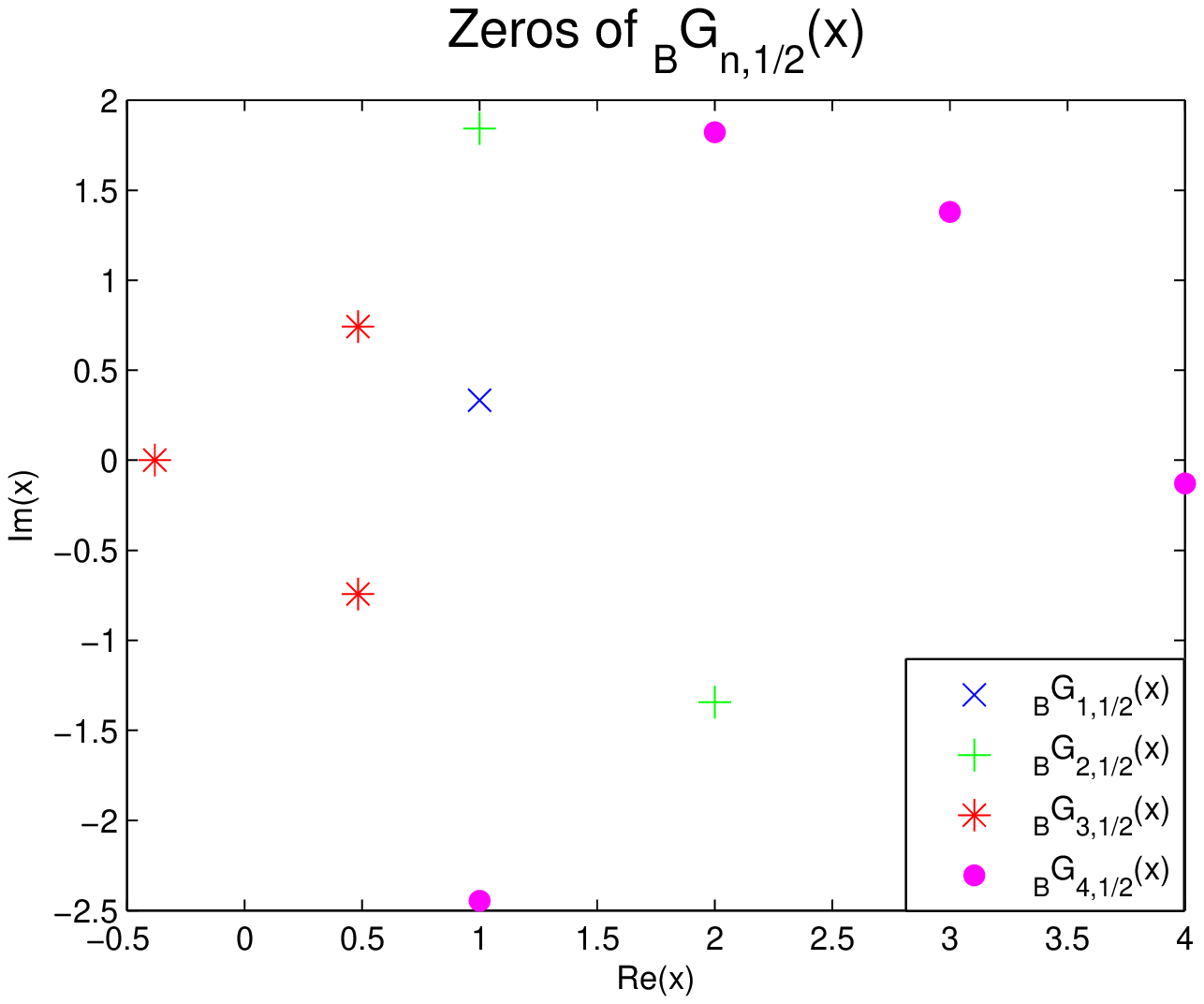}
\includegraphics[width=7.5cm,scale=1.5]{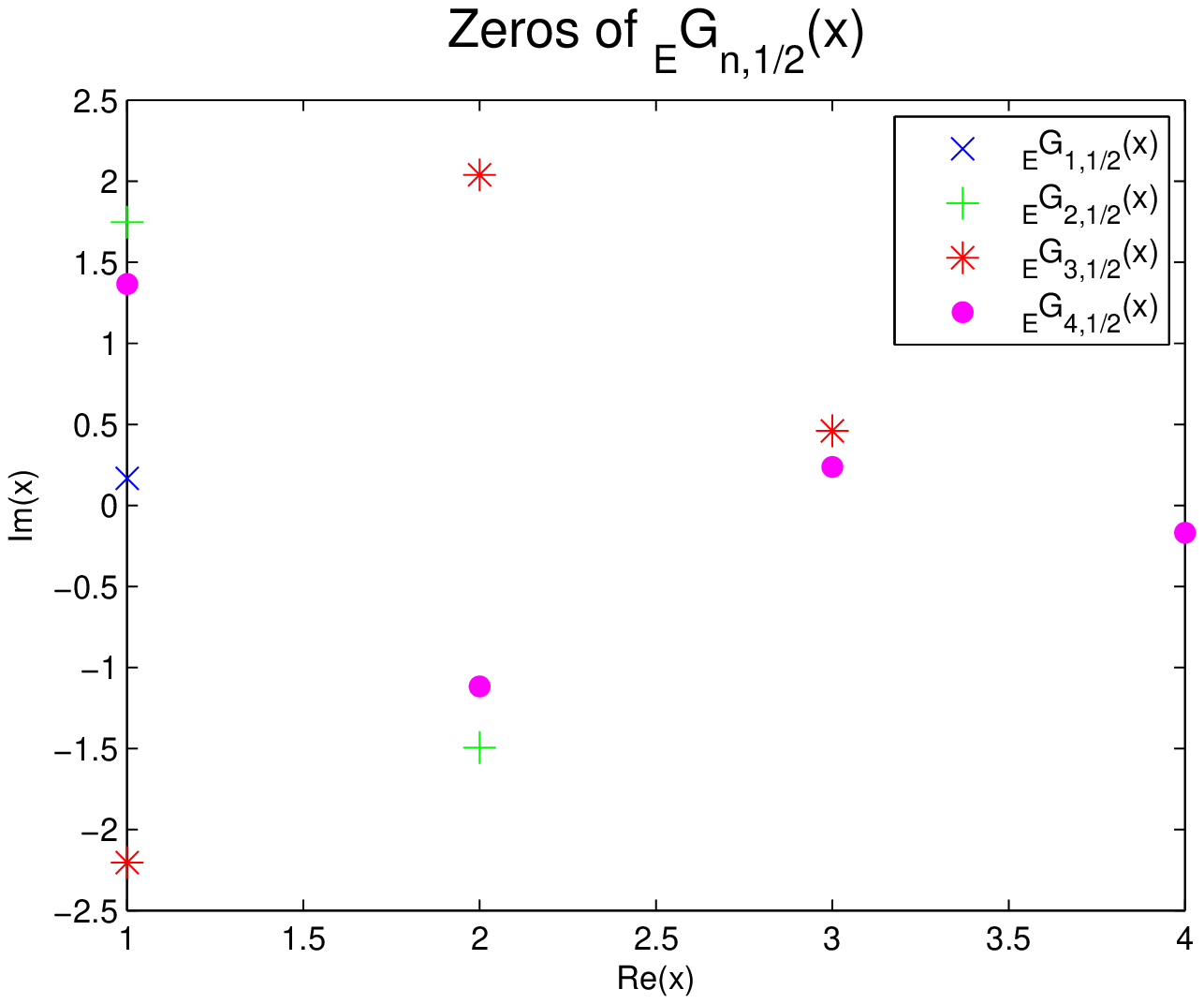}
{\bf Figure 4.11}\hspace{5cm}{\bf Figure 4.12}
\end{center}

\noindent
{\bf Remark~4.1.}~From Tables 9 and 10, the following general relation is observed. The number of real zeros lying on the real plane $\textrm{Im}~(x)=0$, i.e.,
\begin{equation}
\textrm{Real zeros of} ~A_{n,q}^{[2]}(x)=n-\textrm{Complex zeros of} ~A_{n,q}^{[2]}(x),\nonumber
\end{equation}
where $n$ is the degree of polynomial.\\

In order to make the above discussion more clear, we draw the combined graphs of shapes and zeros of the polynomials $B_{n,1/2}^{[2]}(x)$, $E_{n,1/2}^{[2]}(x)$, $G_{n,1/2}^{[2]}(x)$, $_BE_{n,1/2}(x)$, $_BG_{n,1/2}(x)$ and $_EG_{n,1/2}(x)$ for $n=4$.\\

\begin{center}
\includegraphics[width=7.5cm,scale=1.5]{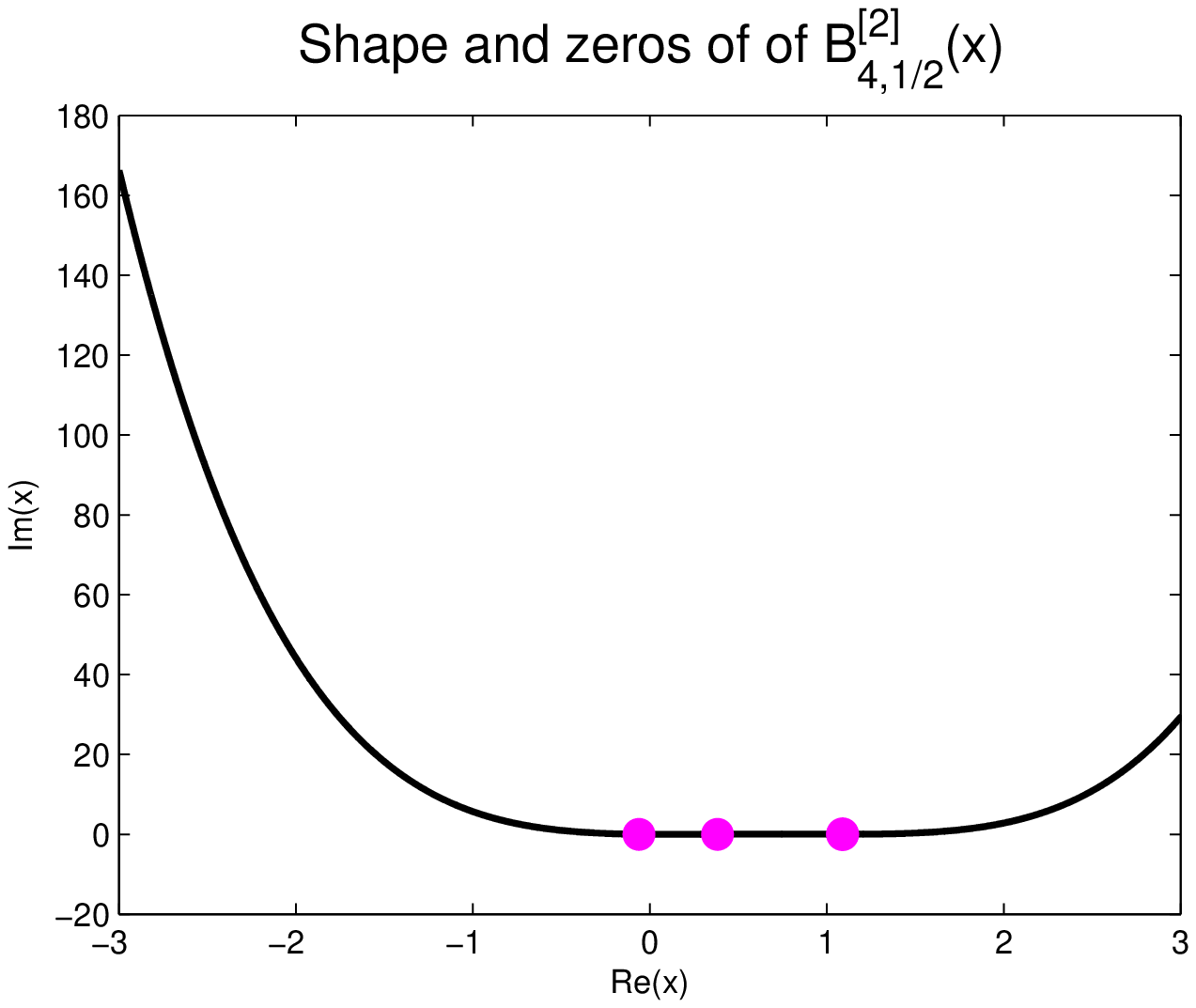}
\includegraphics[width=7.5cm,scale=1.5]{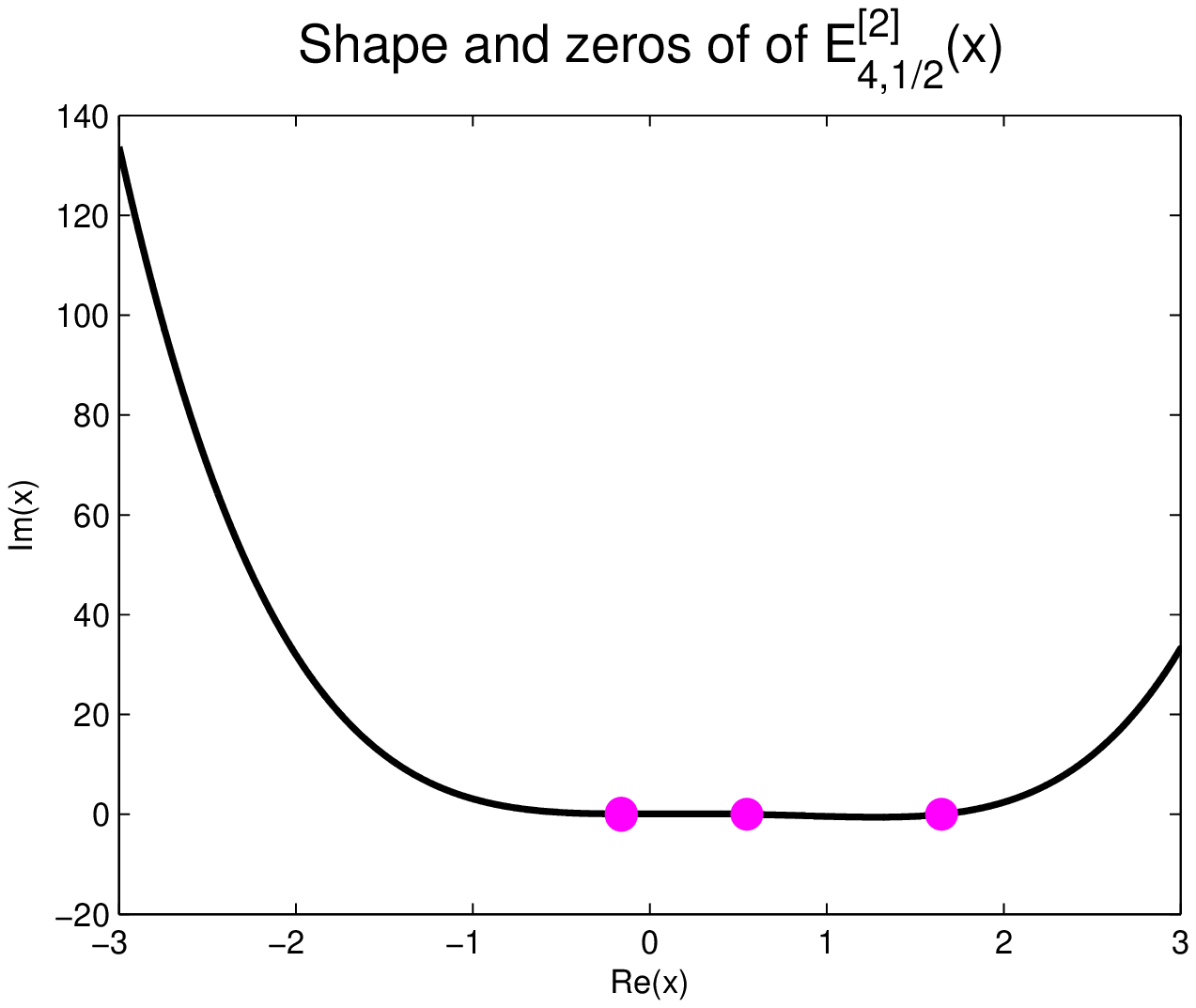}
{\bf Figure 4.13}\hspace{5cm}{\bf Figure 4.14}
\end{center}

\vspace{.5cm}
\begin{center}
\includegraphics[width=7.5cm,scale=1.5]{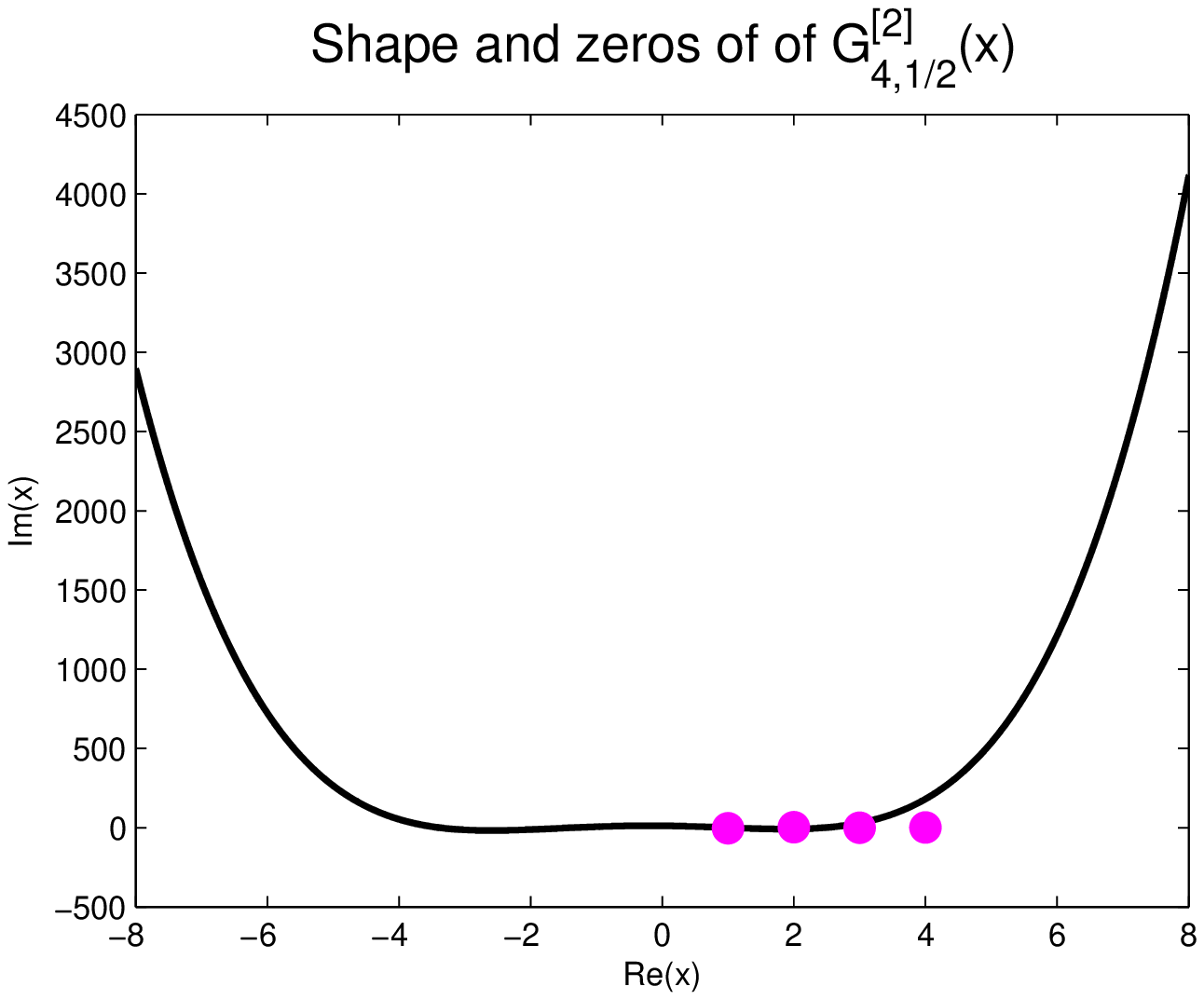}
\includegraphics[width=7.5cm,scale=1.5]{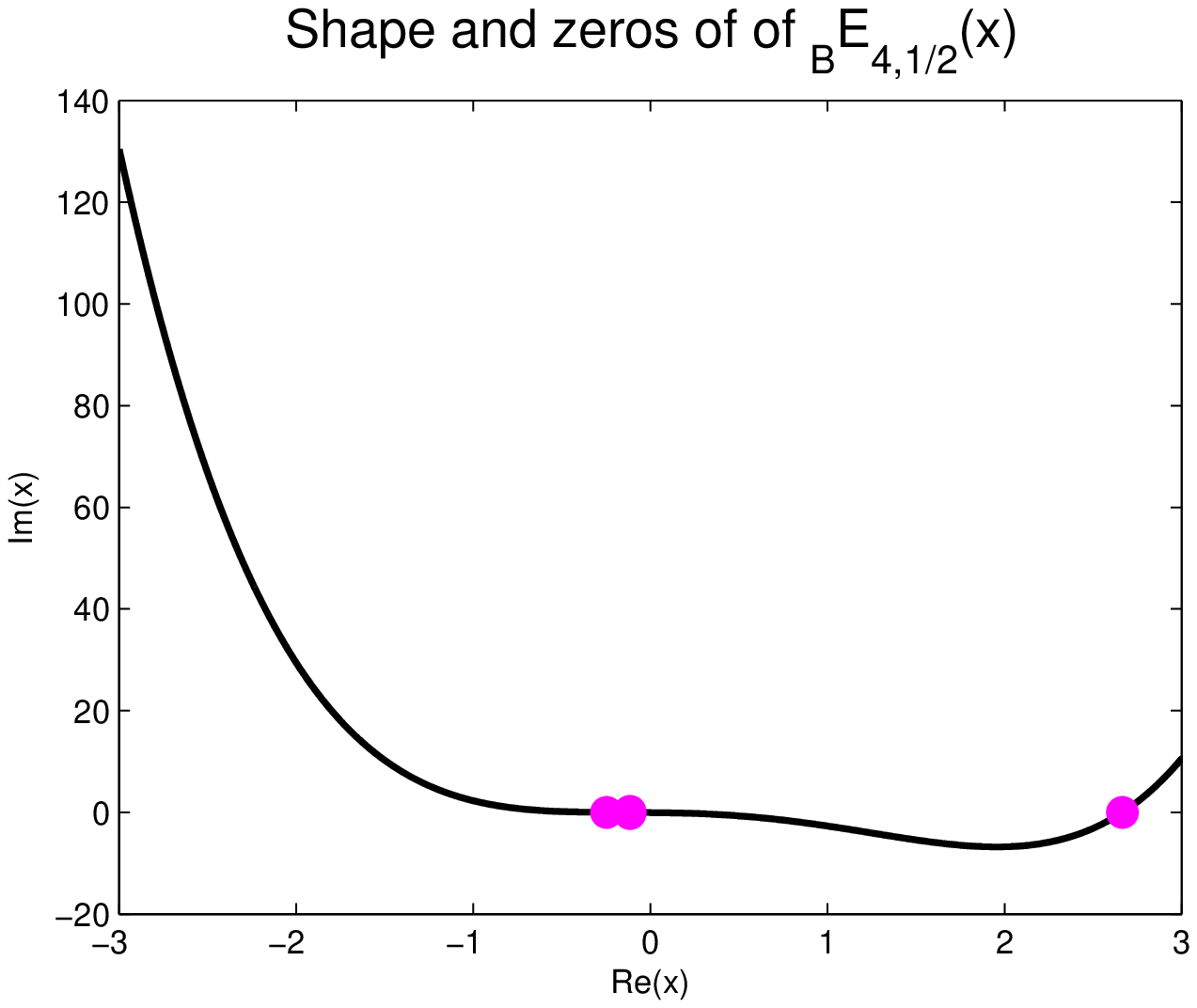}
{\bf Figure 4.15}\hspace{5cm}{\bf Figure 4.16}
\end{center}
\vspace{.5cm}

\begin{center}
\includegraphics[width=7.5cm,scale=1.5]{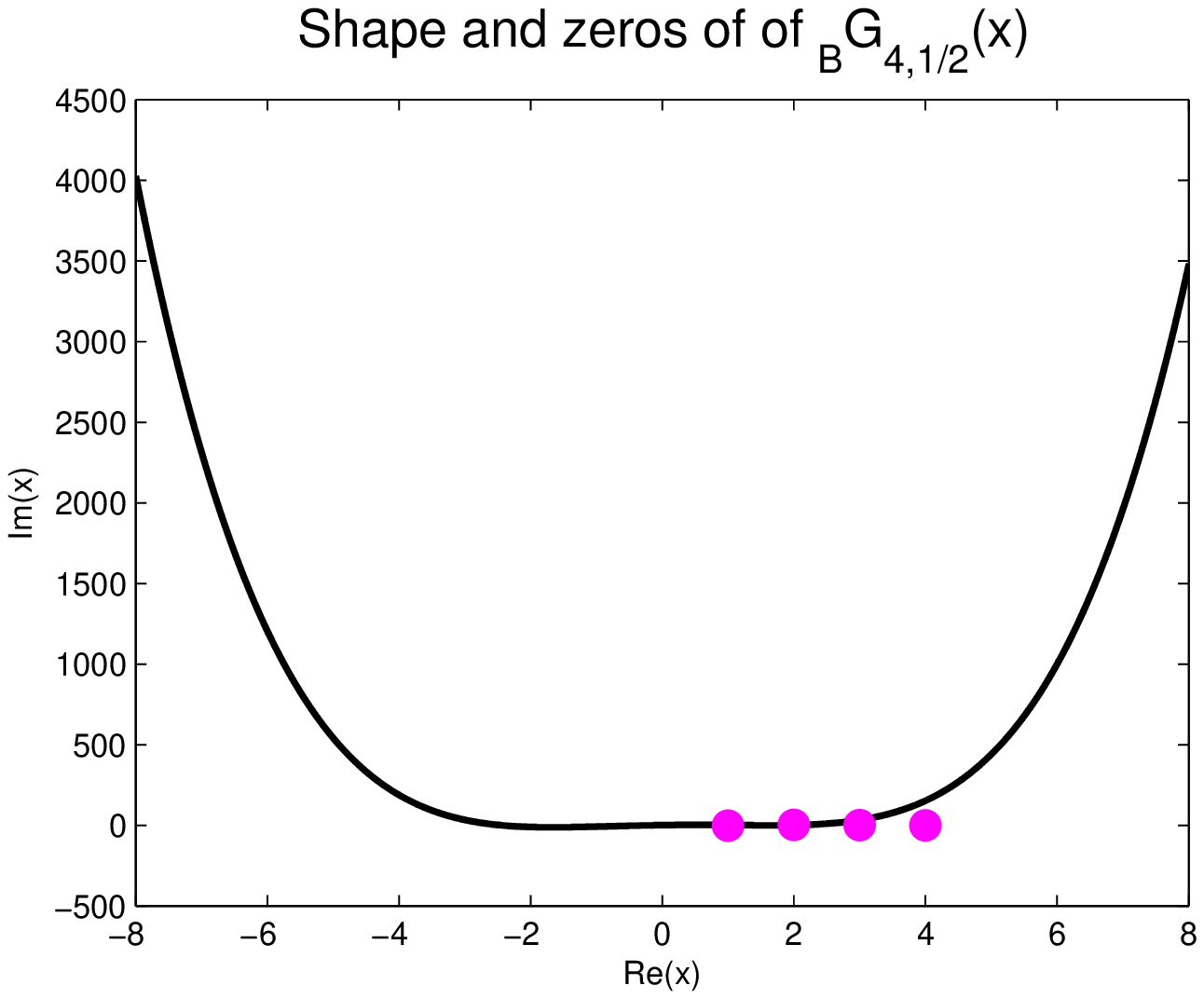}
\includegraphics[width=7.5cm,scale=1.5]{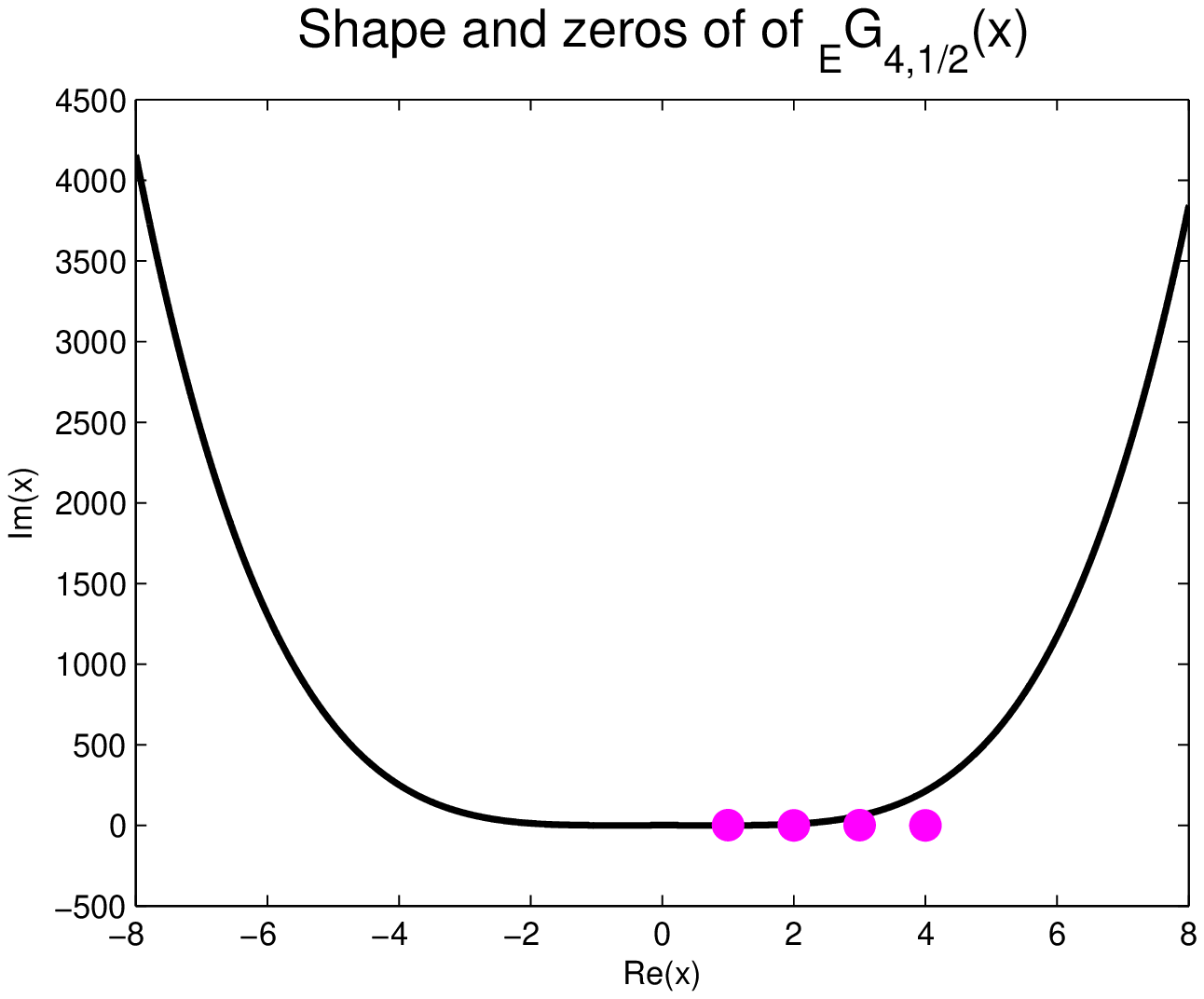}
{\bf Figure 4.17}\hspace{5cm}{\bf Figure 4.18}
\end{center}

It is to be noted that in Figures 4.13, 4.14 and 4.16 out of total two complex zeros only one, with positive imaginary part
is visible, due to the absence of negative imaginary axis in these graphs.\\

The numerical results for approximate solutions of real zeros of $B_{n,1/2}^{[2]}(x)$, $E_{n,1/2}^{[2]}(x)$, $G_{n,1/2}^{[2]}(x)$, $_BE_{n,1/2}(x)$, $_BG_{n,1/2}(x)$ and $_EG_{n,1/2}(x)$ for $(n=1,2,3,4)$ are displayed in Table 9.\\

Also, we note that the real zeros of these polynomials displayed in Table 9 are giving the numerical results for the approximate
solutions of $B_{n,1/2}^{[2]}(x)=0$, $E_{n,1/2}^{[2]}(x)=0$, $G_{n,1/2}^{[2]}(x)=0$, $_BE_{n,1/2}(x)=0$, $_BG_{n,1/2}(x)=0$ and $_EG_{n,1/2}(x)=0$ for $n=1,2,3,4$.\\

\noindent
{\bf Note.}~It is important to note that the 2-iterated $q$-Appell polynomials introduced in this article are actually the $q$-Appell polynomials, since their generating function is of the type $A^{\star}_q(t)e_q(xt)$, with a suitable choice for $A^{\star}_q(t)$. The results established in this article for the 2-iterated and mixed type $q$-special polynomials can be used to solve the existing as well as new emerging problems of certain branches of mathematics, physics and engineering.\\

\noindent
{\large{\bf Acknowledgement}}
\vspace{.35cm}

This work has been done under Senior Research Fellowship (Award letter No. F1-17.1/2012-13, MANF-MUS-UTT-9243) awarded to the second author by the University Grants Commission, Government of India, New Delhi.\\


\begin{thebibliography}{99}
\bibitem{Qappl}W.A. Al-Salam, $q$-Appell polynomials, Ann. Mat. Pura Appl. {\bf 4}(17) (1967) 31-45.
\bibitem{QB}W.A. Al-Salam, $q$-Bernoulli numbers and polynomials, Math. Nachr. {\bf 17} (1959) 239-260.
\bibitem{qMath}G.E. Andrews, R. Askey, R. Roy, 71th Special functions of Encyclopedia of Mathematics and its applications,
Cambridge University Press, Cambridge, 1999.
\bibitem{Appell}P. Appell, Sur une classe de polyn\^{o}mes, Ann. Sci. $\acute{E}$cole. Norm. Sup. {\bf 9}(2) (1880) 119-144.
\bibitem{Ernst1}T. Ernst, $q$-Bernoulli and $q$-Euler polynomials, An umbral approach, Int. J. Diff. Equ. 0973-6069 {\bf 1}(1) (2006) 31-80.
\bibitem{q-Appell}M.E. Keleshteri, N.I. Mahmudov, A study on $q$-Appell polynomials from determinantal point of view, Appl. Math. Comput. {\bf 260} (2015) 351-369.
\bibitem{q-UmbralAppl}M.E. Keleshteri, N.I. Mahmudov, A $q$-umbral approach to $q$-Appell polynomials, arXiv:1505.05067.
\bibitem{qBE} N.I. Mahmudov, On a class of $q$-Bernoulli and $q$-Euler polynomials, Adv. Difference Equ. {\bf 108}
(2013) 11.

\bibitem{Rom1}  S. Roman, The theory of the umbral calculus I, J. Math. Anal. Appl. 87 (1982) 58-115.
\bibitem{Rom} S. Roman, More on the umbral calculus, with emphasis on the $q$-umbral calculus, J. Math. Anal. Appl. {\bf 107} (1985) 222-254.
\bibitem{Ryoo}C.S. Ryoo, A note on $q$-Bernoulli numbers and polynomials, Appl. Math Lett, {\bf 20} (2007) 524-531.
\bibitem{Ryoo1}C.S. Ryoo, T. Kim, R.P. Agarwal, A numerical investigation of the roots of $q$-polynomials, Int. J. Comput. Math. {\bf 83}(2) (2006) 223-234.
\bibitem{SharmaChak} A. Sharma, A.M. Chak, The basic analogue of a class of polynomials, Riv. Mat. Univ. Parma {\bf 5} (1954) 325-337.

\bibitem{SrivqAppl}H. M. Srivastava, Some characterizations of Appell and $q$-Appell polynomials, Ann. Mat. Pura Appl. {\bf 4}(130) (1982) 321-329.
\bibitem{DetSheffer}W. Wang, A determinantal approach to Sheffer sequences, Linear Algebra Appl. {\bf 463} (2014) 228-254.
\end{thebibliography}
\end{document}